# Covering properties of meromorphic functions, negative curvature and spherical geometry

By M. Bonk and A. Eremenko*

## Abstract

Every nonconstant meromorphic function in the plane univalently covers spherical discs of radii arbitrarily close to $\arctan\sqrt{8} \approx 70°32'$. If in addition all critical points of the function are multiple, then a similar statement holds with $\pi/2$. These constants are the best possible. The proof is based on the consideration of negatively curved singular surfaces associated with meromorphic functions.

## 1. Introduction

Let $\mathcal{M}$ be the class of all nonconstant meromorphic functions $f$ in the complex plane $\mathbb{C}$. In this paper we exhibit a universal property of functions $f$ in $\mathcal{M}$ by producing sharp lower bounds for the radii of discs in which branches of the inverse $f^{-1}$ exist. Since a meromorphic function is a mapping into the Riemann sphere $\overline{\mathbb{C}}$, it is appropriate to measure the radii of discs in the spherical metric on $\overline{\mathbb{C}}$. This metric has length element $2|dw|/(1+|w|^2)$ and is induced by the standard embedding of $\overline{\mathbb{C}}$ as the unit sphere $\Sigma$ in $\mathbb{R}^3$. The spherical distance between two points in $\Sigma$ is equal to the angle between the directions to these points from the origin.

Let $D$ be a region in $\mathbb{C}$, and $f\colon D \to \overline{\mathbb{C}}$ a nonconstant meromorphic function. For every $z_0$ in $D$ we define $b_f(z_0)$ to be the spherical radius of the largest open spherical disc centered at $f(z_0)$ for which there exists a holomorphic branch $\phi_{z_0}$ of the inverse $f^{-1}$ with $\phi_{z_0}(f(z_0)) = z_0$. If $z_0$ is a critical point, then $b_f(z_0) := 0$. We define the *spherical Bloch radius* of $f$ by

$$\mathfrak{B}(f) := \sup\{b_f(z_0) : z_0 \in D\},$$

* The first author was supported by a Heisenberg fellowship of the DFG. The second author was supported by NSF grant DMS-9800084 and by Bar-Ilan University.

and the spherical Bloch radius for the class $\mathcal{M}$ by

$$\mathfrak{B} := \inf\{\mathfrak{B}(f) : f \in \mathcal{M}\}.$$

An upper bound for $\mathfrak{B}$ can be obtained from the following example (cf. [23], [20]). We consider a conformal map $f_0$ of an equilateral Euclidean triangle onto an equilateral spherical triangle with angles $2\pi/3$. We always assume that maps between triangles send vertices to vertices. By symmetry $f_0$ has an analytic continuation to a meromorphic function in $\mathbb{C}$. The critical points of $f_0$ form a regular hexagonal lattice and its critical values correspond to the four vertices of a regular tetrahedron inscribed in the sphere $\Sigma$.

If we place one of the vertices of the tetrahedron at the point corresponding to $\infty \in \overline{\mathbb{C}}$ and normalize the map by $z^2 f_0(z) \to 1$ as $z \to 0$, then $f_0$ becomes a Weierstrass $\wp$-function with a hexagonal lattice of periods. It satisfies the differential equation

$$(\wp')^2 = 4(\wp - e_1)(\wp - e_2)(\wp - e_3),$$

where the numbers $e_j$ correspond to the three remaining vertices of the tetrahedron.

It is easy to see that $\mathfrak{B}(f_0) = b_0$, where

$$b_0 := \arctan\sqrt{8} = \arccos(1/3) \approx 1.231 \approx 70^\circ 32'$$

is the spherical circumscribed radius of a spherical equilateral triangle with all angles equal to $2\pi/3$. Hence $\mathfrak{B} \leq b_0$.

Our main result is

THEOREM 1.1. $\mathfrak{B} = b_0$.

The lower estimate $\mathfrak{B} \geq b_0$ in Theorem 1.1 is obtained by letting $R$ tend to infinity in the next theorem. We use the notation $D(R) = \{z \in \mathbb{C} : |z| < R\}$ and $\mathbb{D} = D(1)$.

THEOREM 1.2. *There exists a function* $C_0 : (0, b_0) \to (0, \infty)$ *with the following property. If* $f$ *is a meromorphic function in* $D(R)$ *with* $\mathfrak{B}(f) \leq b_0 - \varepsilon$, *then*

$$\frac{|f'(z)|}{1 + |f(z)|^2} \leq C_0(\varepsilon) \frac{R}{R^2 - |z|^2}. \tag{1}$$

In other words, for every $\varepsilon \in (0, b_0)$ the family of all meromorphic functions on $D(R)$ with the property $\mathfrak{B}(f) \leq b_0 - \varepsilon$ is a normal invariant family [18, 6.4], and each function of this family is a normal function.

The history of this problem begins in 1926 when Bloch [9] extracted the following result from the work of Valiron [26]: *Every nonconstant entire function has holomorphic branches of the inverse in arbitrarily large Euclidean* *discs*. Improving Valiron's arguments he arrived at a stronger statement: *Every holomorphic function $f$ in the unit disc has an inverse branch in some Euclidean disc of radius $\delta|f'(0)|$, where $\delta > 0$ is an absolute constant.* Landau defined Bloch's constant $\mathcal{B}_0$ as the least upper bound of all numbers $\delta$ for which this statement is true. Finding the exact values of $\mathcal{B}_0$ and related constants leads to notoriously hard problems that are mostly unsolved. The latest results for $\mathcal{B}_0$ can be found in [13] and [7]. The conjectured extremal functions for these constants derive from an example given by Ahlfors and Grunsky [6]. As the elliptic function $f_0$ above, the Ahlfors-Grunsky function shows a hexagonal symmetry in its branch point distribution. It seems that our Theorem 1.1 is the first result where a function with hexagonal symmetry is shown to be extremal for a Bloch-type problem.

The earliest estimate for $\mathfrak{B}$ is due to Ahlfors [1], who used what became later known as his Five Islands Theorem (Theorem A below) to prove the lower bound $\mathfrak{B} \geq \pi/4$. We will see that our Theorem 1.1 in turn implies the Five Islands Theorem.

Later Ahlfors [2] introduced another method for treating this type of problems, and obtained a lower bound for Bloch's constant $\mathcal{B}_0$. Applying this method to meromorphic functions, Pommerenke [23] proved an estimate of the form (1) for functions $f$ in $D(R)$ satisfying $\mathfrak{B}(f) \leq \pi/3 - \varepsilon$. From this, one can derive $\mathfrak{B} \geq \pi/3$ thus improving Ahlfors's lower bound. Related is a result by Greene and Wu [16] who showed that for a meromorphic function $f$ in the unit disc the estimate $\mathfrak{B}(f) \leq 18°45'$ implies $|f'(0)|/(1+|f(0)|^2) \leq 1$. An earlier result of this type without numerical estimates is due to Tsuji [25].

Similar problems have been considered for various subclasses of $\mathcal{M}$. In [23] Pommerenke proved that for *locally univalent* meromorphic functions $f$ in $D(R)$ the condition $\mathfrak{B}(f) \leq \pi/2 - \varepsilon$ implies an estimate of the form (1). A different proof was given by Peschl [22]. Minda [19], [20] introduced the classes $\mathcal{M}_m$ of all nonconstant meromorphic functions in $\mathbb{C}$ with the property that all critical points have multiplicity at least $m$. Thus $\mathcal{M}_1 = \mathcal{M}$, $\mathcal{M}_1 \supset \mathcal{M}_2 \supset \ldots \supset \mathcal{M}_\infty$, and $\mathcal{M}_\infty$ is the class of locally univalent meromorphic functions. Using the notation $\mathfrak{B}_m = \inf\{\mathfrak{B}(f) : f \in \mathcal{M}_m\}$, Minda's result can be stated as

$$\mathfrak{B}_m \geq 2\arctan\sqrt{\frac{m}{m+2}}, \quad m \in \mathbb{N} \cup \{\infty\}. \tag{2}$$

In [10] the authors considered some other subclasses. In particular the best possible estimate $\mathfrak{B}(f) \geq \pi/2$ was proved for meromorphic functions omitting at least one value, and $\mathfrak{B}(f) \geq b_0$ was shown for a class of meromorphic functions which includes all elliptic and rational functions.

Since $\mathfrak{B}_1 = \mathfrak{B}$ our Theorem 1.1 improves (2) for $m = 1$. Our method also gives the precise value for all constants $\mathfrak{B}_m$, $m \geq 2$.

Theorem 1.3. *There exists a function $C_1: (0, \pi/2) \to (0, \infty)$ with the following property. If $f$ is a meromorphic function in $D(R)$ with only multiple critical points and $\mathfrak{B}(f) \le \pi/2 - \varepsilon$, then*

$$\frac{|f'(z)|}{1+|f(z)|^2} \le C_1(\varepsilon)\frac{R}{R^2-|z|^2}. \tag{3}$$

*Thus* $\mathfrak{B}_2 = \mathfrak{B}_3 = \ldots = \mathfrak{B}_\infty = \pi/2$.

The first statement of Theorem 1.3 immediately gives the lower bound $\pi/2$ for $\mathfrak{B}_2, \ldots, \mathfrak{B}_\infty$. This bound is achieved as the exponential function $\exp \in \mathcal{M}_\infty$ shows.

The Ahlfors Five Islands Theorem is

Theorem A. *Given five Jordan regions on the Riemann sphere with disjoint closures, every nonconstant meromorphic function $f: \mathbb{C} \to \overline{\mathbb{C}}$ has a holomorphic branch of the inverse in one of these regions.*

*Derivation of Theorem* A *from Theorem* 1.1. We consider the following five points on the Riemann sphere

$$e_1 = \infty,\ e_2 = 0,\ e_3 = 1, \text{ and } e_{4,5} = \exp(\pm 2\pi i/3).$$

These points serve as vertices of a triangulation of the sphere into six spherical triangles, each having angles $\pi/2$, $\pi/2$, $2\pi/3$. The spherical circumscribed radius of each of these triangles is $R_0 := \arctan 2 \approx 63°26'$. (See for example [12, p. 246].) This means that each point on the sphere is within distance $R_0$ from one of the points $e_j$. Let $\psi: \overline{\mathbb{C}} \to \overline{\mathbb{C}}$ be a diffeomorphism which sends the given Jordan regions $D_j$, $1 \le j \le 5$, into the spherical discs $B_j$ of radius $\varepsilon_0 := (b_0 - R_0)/2 > 0$ centered at $e_j$, $1 \le j \le 5$. By the Uniformization Theorem there exists a quasiconformal diffeomorphism $\phi: \mathbb{C} \to \mathbb{C}$ and a meromorphic function $g: \mathbb{C} \to \overline{\mathbb{C}}$ such that

$$\psi \circ f = g \circ \phi. \tag{4}$$

By Theorem 1.1 an inverse branch of $g$ exists in some spherical disc $B$ of radius $b_0 - \varepsilon_0$. Every such disc $B$ contains at least one of the discs $B_j$. So $g$ and thus $\psi \circ f$ have inverse branches in one of the discs $B_j$. We conclude that $f$ has an inverse branch in one of the regions $D_j \subset \psi^{-1}(B_j)$. □

The use of the diffeomorphism $\psi$ in the above proof was suggested by recent work of Bergweiler [8] who gives a simple proof of the Five Islands Theorem using a normality argument.

Our Theorem 1.2 implies a stronger version of the Five Islands Theorem proved by Dufresnoy [14], [18].

THEOREM B. *Let $D_1, \ldots, D_5$ be five Jordan regions on the Riemann sphere whose closures are disjoint. Then there exists a positive constant $C_2$, depending on these regions, with the following property. Every meromorphic function $f$ in $\mathbb{D}$ without inverse branches in any of the regions $D_j$ satisfies*

$$\frac{|f'(0)|}{1+|f(0)|^2} \leq C_2.$$

In fact one can take $C_2 = (32L\max\{C_0(\varepsilon_0), 1\})^K$, where $L \geq 1$ is the Lipschitz constant of $\psi^{-1}$ in (4), $K$ is the maximal quasiconformal dilatation of $\psi$, $\varepsilon_0$ is as above and $C_0$ is the function from Theorem 1.2. So $C_0(\varepsilon_0)$ is an absolute constant, while $L$ and $K$ depend on the choice of the regions $D_j$ in Theorem B. This value for $C_2$ can be obtained by an application of Mori's theorem similarly as below in our reduction of Theorems 1.2 and 1.3 to Theorem 1.4.

Before we begin discussing the proofs of Theorems 1.2 and 1.3, let us introduce some notation and fix our terminology. It is convenient to use the language of singular surfaces though our surfaces are of very simple kind, called $K$-polyhedra in [24, Ch. I, 5.7]. For $r \in (0,1]$, $\alpha > 0$ and $\chi \in \{0,1,-1\}$ a *cone* $C(\alpha,\chi,r)$, is the disc $D(r)$, equipped with the metric given by the length element

$$\frac{2\alpha|z|^{\alpha-1}|dz|}{1+\chi|z|^{2\alpha}}.$$

This metric has constant Gaussian curvature $\chi$ in $D(r)\backslash\{0\}$.

To visualize a cone we choose a sequence $0 = \alpha_0 < \alpha_1 < \ldots < \alpha_n = 2\pi\alpha$ with $\alpha_j - \alpha_{j-1} < 2\pi$, and $r \in (0,1]$, and consider the closed sectors

$$D_j = \{w \in D(r^\alpha) : \alpha_{j-1} \leq \arg w \leq \alpha_j\},\ 1 \leq j \leq n,$$

equipped with the Riemannian metric of constant curvature $\chi$, whose length element is $2|dw|/(1+\chi|w|^2)$. For $1 \leq j \leq n-1$ we paste $D_j$ to $D_{j+1}$ along their common side $\{w \in D(r^\alpha) : \arg w = \alpha_j\}$, and then identify the remaining side $\{w \in D(r^\alpha) : \arg w = \alpha_n\}$ of $D_n$ with the side $\{w \in D(r^\alpha) : \arg w = \alpha_0\}$ of $D_1$, all identifications respecting arclength. Thus we obtain a singular surface $S$ which is isometric to the cone $C(\alpha,\chi,r)$ via $z = \phi(w) = w^{1/\alpha}$.

We consider a two-dimensional connected oriented triangulable manifold (a surface) equipped with an intrinsic metric, which means that the distance between every two points is equal to the infimum of lengths of curves connecting these points. By a *singular surface* we mean in this paper a surface with an intrinsic metric which satisfies the following condition. For every point $p$ there exists a neighborhood $V$ and an isometry $\phi$ of $V$ onto a cone $C(\alpha,\chi,r)$. The numbers $r$ and $\alpha$ in the definition of a cone may vary from one point to another. It follows from this definition that near every point, except some isolated set of *singularities* we have a smooth Riemannian metric of constant curvature $\chi \in \{0, 1, -1\}$. The curvature at a singular point $p$ is defined to be $+\infty$ if $0 < \alpha < 1$ and $-\infty$ if $\alpha > 1$. The *total angle* at $p$ is $2\pi\alpha$, and $p$ contributes $2\pi(1 - \alpha)$ to the integral curvature.

Underlying the metric structure of a singular surface is a canonical Riemann surface structure. It is obtained by using the local coordinates $\phi$ from the definition of a singular surface as conformal coordinates. When we speak of a 'conformal map' or 'uniformization' of singular surfaces, we mean this conformal structure. So our 'conformal maps' do not necessarily preserve angles at singularities.

In what follows every hyperbolic region in the plane is assumed to carry its unique smooth complete Riemannian metric of constant curvature $-1$, unless we equip it explicitly with some other metric. For example, $D(R)$ is always assumed to have the metric with length element

$$\frac{2R|dz|}{R^2 - |z|^2}.$$

If $\phi: D(R) \to Y$ is a conformal map of singular surfaces, and the curvature on $Y$ is at most $-1$ everywhere, then $\phi$ is distance decreasing. This follows from Ahlfors's extension of Schwarz's lemma [5, Theorem 1-7].

If $f: X \to Y$ is a smooth map between singular surfaces, we will denote by $\|f'\|$ the norm of the derivative with respect to the metrics on $X$ and $Y$. So (1), for example, can simply be written as $\|f'\| \le C_0(\varepsilon)$. We reserve the notation $|f'|$ for the case when the Euclidean metric is considered in both $X$ and $Y$.

In this paper *triangle* always means a triangle whose angles are strictly between 0 and $\pi$, and *spherical triangle* refers to a triangle isometric to one on the unit sphere $\Sigma$ in $\mathbb{R}^3$.

Let $D$ be a region in $\overline{\mathbb{C}}$, and $f: D \to \overline{\mathbb{C}}$ a nonconstant meromorphic function. We consider another copy of $D$ and convert it into a new singular surface $S_f$ by providing it with the pullback of the spherical metric via $f$. Then the metric on $S_f$ has the length element $2|f'(z)dz|/(1 + |f(z)|^2)$. The identity map id: $D \to D$ now becomes a conformal homeomorphism $f_1: D \to S_f$. Thus $f$ factors as

$$f = f_2 \circ f_1, \quad D \xrightarrow{f_1} S_f \xrightarrow{f_2} \overline{\mathbb{C}}, \tag{5}$$

where $f_2$ is a path isometry, that is, $f_2$ preserves the arclength of every rectifiable path. In the classical literature, the singular surface $S_f$ is called *the Riemann surface of* $f^{-1}$ *spread over the sphere* (Überlagerungsfläche).

The following idea was first used by Ahlfors in his paper [2] to obtain a lower bound for the classical Bloch constant $\mathcal{B}_0$. In the papers of Pommerenke [23] and Minda [19], [20] essentially the same method was applied to the case of the spherical metric. Assuming that the Bloch radius of a function $f: D(R) \to \overline{\mathbb{C}}$ is small enough one constructs a conformal metric on $S_f$ whose curvature is bounded from above by a negative constant. If we denote the surface $S_f$ equipped with this new metric by $S'$, then the identity map $\mathrm{id}: S_f \to S_f$ can be considered as a conformal map $\psi: S_f \to S'$. If, in addition, one knows a lower bound for the norm of the derivative of $\psi$, an application of the Ahlfors-Schwarz lemma to $\psi \circ f_1$ (where $f_1$ is as in (5)) leads to an estimate of the form (1). For a proof of Theorem 1.2 it seems hard to find an explicit conformal map $\psi$ which will work in the case when $\mathfrak{B}(f)$ is close to $b_0$.

Our main innovation is replacing the conformal map $\psi$ by a quasiconformal map such that $\psi^{-1} \in \mathrm{Lip}(L)$, for some $L \geq 1$, which means that $\psi^{-1}$ satisfies a Lipschitz condition with constant $L$ (see Theorem 1.4 below).

Let us describe the scheme of our proof.

We recall that $f: D \to \overline{\mathbb{C}}$ is said to have an asymptotic value $a \in \overline{\mathbb{C}}$ if there exists a curve $\gamma: [0, 1) \to D$ such that $\gamma(t)$ leaves every compact subset of $D$ as $t \to 1$ and $a = \lim_{t\to 1} f(\gamma(t))$. The potential presence of asymptotic values causes difficulties, so our first step is a reduction of Theorems 1.2 and 1.3 to their special cases for functions without asymptotic values. This reduction is based on a simple approximation argument (Lemmas 2.1 and 2.2). If $f$ has no asymptotic values, then the singular surface $S_f$ is complete; that is, every curve of finite length in $S_f$ has a limit in $S_f$.

As a second step we introduce a locally finite covering[1] $T$ of $S_f$ by closed spherical triangles such that the intersection of any two triangles of $T$ is either empty or a common side or a set of common vertices. In addition the set of vertices of these triangles coincides with the critical set of $f_2$, and the circumscribed radii of all triangles do not exceed the Bloch radius $\mathfrak{B}(f)$. The existence of such covering was proved in [10] under the conditions that $f$ has no asymptotic values, and $\mathfrak{B}(f) < \pi/2$. For the precise formulation see Lemma 2.3 in Section 2. Our singular surface $S_f$ has spherical geometry everywhere, except at the vertices of the triangles where it has singularities. Since each vertex of a triangle in $T$ is a critical point of $f_2$ in (5), the total angle at a vertex $p$ is $2\pi m$ where $m$ is the local degree of $f_2$ at $p$, $m \geq 2$. In particular, the total angle at each vertex of a triangle in $T$ is at least $4\pi$, and at least $6\pi$ if all critical points of $f$ are multiple. Now our results will follow from

THEOREM 1.4. *Let $S$ be an open simply-connected complete singular surface with a locally finite covering by closed spherical triangles such that the intersection of any two triangles is either empty, a common side, or a set of*

[1] This covering need not be a triangulation, since two distinct triangles might have more than one common vertex without sharing a common side.

*common vertices. Assume that for some* $\varepsilon > 0$ *one of the following conditions holds*:

(i) *The circumscribed radius of each triangle is at most* $b_0 - \varepsilon$ *and the total angle at each vertex is at least* $4\pi$, *or*

(ii) *The circumscribed radius of each triangle is at most* $\pi/2 - \varepsilon$ *and the total angle at each vertex is at least* $6\pi$.

*Then there exists a* $K$*-quasiconformal map* $\psi: S \to \mathbb{D}$ *such that* $\psi^{-1} \in \mathrm{Lip}(L)$, *with* $L$ *and* $K$ *depending only on* $\varepsilon$.

Case (i) will give Theorem 1.2 and case (ii) will give Theorem 1.3.

The idea behind Theorem 1.4 is that if the triangles are small enough, then the negative curvature concentrated at vertices dominates the positive curvature spread over the triangles. Thus on a large scale our surface looks like one whose curvature is bounded from above by a negative constant.

For the proof in the case (i) we first construct a new singular surface $\tilde{S}$ in the following way. We choose an appropriate increasing subadditive function $F: [0, \pi) \to [0, \infty)$ with $F(0) = 0$ and

$$F'(0) < \infty, \tag{6}$$

and replace each spherical triangle $\Delta \in T$ with sides $a, b, c$ by a Euclidean triangle $\tilde{\Delta}$ whose sides are $F(a), F(b), F(c)$. The monotonicity and subadditivity of $F$ imply that this is possible for each triangle $\Delta$; that is, $F(a), F(b), F(c)$ satisfy the triangle inequality. Then these new triangles $\tilde{\Delta}$ are pasted together according to the same combinatorial pattern as the triangles $\Delta$ in $T$, with identification of sides respecting arclength. Thus we obtain a new singular surface $\tilde{S}$, and (6) permits us to define a bilipschitz homeomorphism $\psi_1: S \to \tilde{S}$.

The main point is to choose $F$ so that each angle of every triangle $\tilde{\Delta}$ is at least $1/2 + \delta$ times the corresponding angle of $\Delta$, where $\delta > 0$ is a constant depending only on $\varepsilon$ (cf. Theorem 1.5). Thus the new singular surface $\tilde{S}$ has Euclidean geometry everywhere, except at the vertices of the triangles, where the total angle is at least $2\pi + 4\pi\delta$. As the diameters of the triangles $\tilde{\Delta}$ are bounded, it is relatively easy (using the Ahlfors method described above) to show that $\tilde{S}$ is conformally equivalent to $\mathbb{D}$, and that the uniformizing conformal map $\psi_2: \tilde{S} \to \mathbb{D}$ has an inverse in $\mathrm{Lip}(L)$ with $L$ depending only on $\delta$. Then Theorem 1.4 with (i) follows with $\psi = \psi_2 \circ \psi_1$. Case (ii) is treated similarly.

The major difficulty is to verify that some function $F$ has all the necessary properties. To formulate our main technical result we use the following notation. Let $\Delta$ be a spherical or Euclidean triangle whose sides have lengths $a, b, c$, and $F$ a subadditive increasing function with $F(0) = 0$ and $F(t) > 0$ for $t>0$. We define the transformed triangle $\tilde{\Delta} := F\Delta$ as a Euclidean triangle with sides $F(a), F(b), F(c)$. Let $\alpha,\beta,\gamma$ be the angles of $\Delta$, and $\tilde{\alpha},\tilde{\beta},\tilde{\gamma}$ the corresponding angles of $\tilde{\Delta}$. We define the angle distortion of $\Delta$ under $F$ by

$$D(F,\Delta)=\min\{\tilde{\alpha}/\alpha,\tilde{\beta}/\beta,\tilde{\gamma}/\gamma\}. \tag{7}$$

For case (i) in Theorem 1.4 we take $k>0$ and put[2]

$$F_k(t)=\min\{k\,\mathrm{chd}\,t,\sqrt{\mathrm{chd}\,t}\},\quad \text{where}\quad \mathrm{chd}\,t := 2\sin(t/2),\quad t\in[0,\pi]. \tag{8}$$

THEOREM 1.5. *For every $\varepsilon\in(0,b_0)$ there exist $k\geq 1$ and $\delta>0$ such that for every spherical triangle $\Delta$ of circumscribed radius at most $b_0-\varepsilon$ the angle distortion by $F_k$ satisfies*

$$D(F_k,\Delta)\geq \frac{1}{2}+\delta.$$

For Theorems 1.3 and 1.4 with condition (ii) we need a simpler result with

$$F_k^*(t)=\min\{k\,\mathrm{chd}\,t,1\},\quad t\in[0,\pi]. \tag{9}$$

THEOREM 1.6. *For every $\varepsilon\in(0,\pi/2)$ there exist $k\geq 1$ and $\delta>0$ such that for every spherical triangle $\Delta$ of circumscribed radius at most $\pi/2-\varepsilon$ the angle distortion by $F_k^*$ satisfies*

$$D(F_k^*,\Delta)\geq \frac{1}{3}+\delta.$$

The constant $1/2$ in Theorem 1.5 is the best possible, no matter which distortion function $F$ is applied to the sides. Indeed, a spherical equilateral triangle of circumscribed radius close to $b_0$ has a sum of angles close to $2\pi$, and the corresponding Euclidean triangle has a sum of angles $\pi$. A similar remark applies to the constant $1/3$ in Theorem 1.6. Further comments about Theorems 1.5 and 1.6 are in the beginning of Section 3.

*Remarks.* 1. Our proofs of Theorems 1.5 and 1.6 are similar but separate. It seems natural to conjecture that one can 'interpolate' somehow between these results. This would yield Theorem 1.4 under the following condition:

(iii) *For fixed $q\in(1,3]$ and $\varepsilon>0$ the total angle at each vertex is at least $2\pi q$ and the supremum of the circumscribed radii is at most*

$$\arctan\sqrt{\frac{-\cos(\pi q/2)}{\cos^3(\pi q/6)}}-\varepsilon.$$

[2] We find the notation chd (abbreviation of 'chord') convenient. According to van der Waerden [27] the ancient Greeks used the chord as their only trigonometric function. Only in the fifth century were the sine and other modern trigonometric functions introduced.

If $\varepsilon = 0$, then this expression is the circumscribed radius of an equilateral spherical triangle with angles $\pi q/3$. In Theorem 1.4 case (i) corresponds to $q = 2$ and (ii) to $q = 3$. The limiting case $q \to 1$ is also understood, namely we have the

PROPOSITION 1.7. *Every open simply-connected singular surface triangulated into Euclidean triangles of bounded circumscribed radius, and having total angle at least $2\pi q$ at each vertex, where $q > 1$, is conformally equivalent to the unit disc.*

We will prove this in Section 2 as a part of our derivation of Theorem 1.4.

2. Ahlfors's original proof of Theorem A [3], [18] was based on a linear isoperimetric inequality. Assuming that $f$ has no inverse branches in any of the five given regions $D_j$, he deduced that the surface $S_f$ has the property that each Jordan region in $S_f$ of area $A$ and boundary length $L$ satisfies $A \leq hL$, where $h$ is a positive constant depending only on the regions $D_j$. On the other hand, Ahlfors showed that a surface with such a linear isoperimetric inequality cannot be conformally equivalent to the plane, and used this contradiction to prove Theorem A. The stronger Theorem B can be derived by improving the above isoperimetric inequality to

$$A \leq \min\{h_1 L^2, hL\},$$

where the constants $h_1$ and $h$ still depend only on the regions $D_j$. This argument belongs to Dufresnoy [14]; see also [18, Ch. 6]. It seems that a linear isoperimetric inequality holds under any of the conditions of Theorem 1.4 or of the above proposition. This would imply that the surface is hyperbolic in the sense of Gromov [17].

3. It can be shown [10, Lemma 7.2] that condition (i) in Theorem 1.4 implies that the areas of the triangles are bounded away from $\pi$. Similarly (ii) implies that these areas are bounded away from $2\pi$. One can be tempted to replace our assumptions on the circumscribed radii in Theorem 1.4 by weaker (and simpler) assumptions on the areas of the triangles. We believe that there are counterexamples to such stronger versions of Theorem 1.4.

4. Theorem 1.2 can be obtained from Theorem 1.1 by a general rescaling argument as in [10]. (A similar argument derives Bloch's theorem from Valiron's theorem and Theorem B from Theorem A. See [28] for a general discussion of such rescaling.) But we could not simplify our proof by proving the weaker Theorem 1.1 first. The problem is in the crucial approximation argument in Section 2 which deals with asymptotic values.

5. The functions $C_0$ in Theorem 1.2 and $C_1$ in Theorem 1.3 can easily be expressed in terms of $k = k(\varepsilon)$ and $\delta = \delta(\varepsilon)$ from Theorems 1.5 and 1.6. The authors believe that the arguments of this paper can be extended to give explicit estimates for $k$ and $\delta$, but this would make the proofs substantially longer. So we use 'proof by contradiction' to simplify our exposition.

The plan of the paper is the following. In Section 2 we reduce all our results to Theorems 1.5 and 1.6. Section 3 begins with a discussion and outline of the proof of Theorem 1.5. The proof itself occupies the rest of Section 3 as well as Sections 4, 5 and 6. In Section 7 we prove Theorem 1.6, using some lemmas from Sections 4 and 5. The main results of this paper have been announced in [11].

The authors thank D. Drasin, A. Gabrielov and A. Weitsman for helpful discussions, and the referee for carefully reading the paper. A. Gabrielov suggested the use of convexity to simplify the original proof of Lemma 5.1. This paper was written while the first author was visiting Purdue University. He thanks the faculty and staff for their hospitality.

## 2. Derivation of Theorems 1.2, 1.3 and 1.4 from Theorems 1.5 and 1.6

We begin with the argument which permits us to reduce our considerations to functions without asymptotic values.

LEMMA 2.1. *Let $R > 1$, $\varepsilon > 0$, and a meromorphic function $f\colon D(R) \to \overline{\mathbb{C}}$ be given. Then there exists a conformal map $\phi$ of $\mathbb{D}$ into $D(R)$ with $\phi(0) = 0$ and $|\phi'(0) - 1| < \varepsilon$ such that $f \circ \phi$ is the restriction of a rational function to $\mathbb{D}$.*

*If all critical points of $f$ in $\overline{\mathbb{D}}$ are multiple, then $\phi$ can be chosen so that all critical points of $f \circ \phi$ in $\mathbb{D}$ are multiple.*

*Proof.* We may assume that $f$ is nonconstant. Then we can choose an annulus $A = \{z : r_1 < |z| < r_2\}$ with $1 < r_1 < r_2 < R$ such that $f$ has no poles and no critical points in $\overline{A}$. Put $m = \min\{|f'(z)|\colon z \in \overline{A}\} > 0$, and let $\delta \in (0, m/6)$. Taking a partial sum of the Laurent series of $f$ in $A$, we obtain a rational function $g$ such that

$$(10) \qquad |g(z) - f(z)| < (r_2 - r_1)\delta \quad \text{and} \quad |g'(z) - f'(z)| < \delta \quad \text{for} \quad z \in \overline{A}.$$

Let $\lambda$ be a smooth function defined on a neighborhood of $[r_1, r_2]$ such that $0 \le \lambda \le 1$, $\lambda(r) = 0$ for $r \le r_1$, $\lambda(r) = 1$ for $r \ge r_2$, and $|\lambda'| \le 2/(r_2 - r_1)$. We put $u(z) = \lambda(|z|)(g(z) - f(z))$ for $z$ in a neighborhood of $\overline{A}$. Then (10) implies $|u'(z)| \le 3\delta \le m/2$ for $z \in \overline{A}$. The function

$$h(z) = \begin{cases} f(z), & |z| < r_1, \\ f(z) + u(z), & r_1 \le |z| \le r_2, \\ g(z), & |z| > r_2 \end{cases}$$

is meromorphic in $\overline{\mathbb{C}} \backslash \overline{A}$. In a neighborhood of $\overline{A}$ it is smooth, and for its Beltrami coefficient $\mu_h = h_{\overline{z}}/h_z$ we obtain $|\mu_h| \le 6\delta/m < 1$. Thus $h: \overline{\mathbb{C}} \to \overline{\mathbb{C}}$ is a quasiregular map. Hence there exists a quasiconformal homeomorphism $\phi: \overline{\mathbb{C}} \to \overline{\mathbb{C}}$ fixing $0, 1$ and $\infty$ with $\mu_{\phi^{-1}} = \mu_h$. Then $h \circ \phi$ is a rational function. Moreover, when $\delta$ is small, then $\phi$ is close to the identity on $\overline{\mathbb{C}}$. Hence for sufficiently small $\delta > 0$ the homeomorphism $\phi$ is conformal in $\mathbb{D}$, and we have $|\phi'(0) - 1| \le \varepsilon$ and $\phi(\mathbb{D}) \subseteq D(r_1)$. Moreover, if $f$ has only multiple critical points in $\overline{\mathbb{D}}$, then we may in addition assume that $\phi$ maps $\mathbb{D}$ into a disc in which $f$ has only multiple critical points. Thus $f \circ \phi$ has only multiple critical points in $\mathbb{D}$. □

LEMMA 2.2. *Given $\varepsilon > 0$ and a meromorphic function $f: \mathbb{D} \to \overline{\mathbb{C}}$ there exists a meromorphic function $g: \mathbb{D} \to \overline{\mathbb{C}}$ without asymptotic values such that*

$$\mathfrak{B}(g) \le \mathfrak{B}(f) + \varepsilon \tag{11}$$

*and*

$$\|g'(0)\| \ge (1 - \varepsilon)\|f'(0)\|. \tag{12}$$

*If all critical points of $f$ are multiple, then $g$ can be chosen so that all its critical points are also multiple.*

*Proof.* First we approximate $f$ by a restriction of a rational function to $\mathbb{D}$. Assuming $0 < \varepsilon < 1/2$ we set $r = 1 - \varepsilon/2$ and apply the previous Lemma 2.1 to the function $f_r(z) := f(rz)$, $z \in \mathbb{D}$, which is meromorphic in $D(1/r)$. We obtain a conformal map $\phi$ on $\mathbb{D}$ with the properties $\phi(0) = 0$, $\phi(\mathbb{D}) \subset D(1/r)$, and

$$|\phi'(0)| \ge 1 - \varepsilon/2$$

such that $p := f_r \circ \phi: \mathbb{D} \to \overline{\mathbb{C}}$ is the restriction of a rational function $h: \overline{\mathbb{C}} \to \overline{\mathbb{C}}$. (The reason why we have to distinguish between $p$ and $h$ is that in general $\mathfrak{B}(p) \ne \mathfrak{B}(h)$.) We have

$$\mathfrak{B}(p) \le \mathfrak{B}(f) \tag{13}$$

and

$$\|h'(0)\| = \|p'(0)\| \ge (1 - \varepsilon)\|f'(0)\|. \tag{14}$$

If all critical points of $f$ in $\mathbb{D}$ are multiple, then $f_r$ has only multiple critical points in $\overline{\mathbb{D}}$. Then Lemma 2.1 ensures that $\phi$ can be chosen in such a way that $p$ has only multiple critical points in $\mathbb{D}$.

Now we will replace $p$ by a function $g: \mathbb{D} \to \overline{\mathbb{C}}$ which has no asymptotic values.

We consider the singular surface $S_h$ obtained by providing $\overline{\mathbb{C}}$ with the pullback of the spherical metric via $h$. Then $h$ factors as in (5), namely $h = h_2 \circ h_1$, where $h_1: \overline{\mathbb{C}} \to S_h$ is the natural homeomorphism and $h_2: S_h \to \overline{\mathbb{C}}$ is a path isometry.

The compact set $K := S_h \backslash h_1(\mathbb{D})$ has a finite $\varepsilon/2$-net $E \subset K$; that is, every point of $K$ is within distance of $\varepsilon/2$ from $E$. We may assume without loss of generality that $E$ contains at least 4 points. We put

$$H := h_1^{-1}(E) \cup (\mathrm{crit}(h) \cap (\overline{\mathbb{C}} \backslash \mathbb{D})) \subset \overline{\mathbb{C}} \backslash \mathbb{D}, \tag{15}$$

where $\mathrm{crit}(h)$ stands for the set of critical points of $h$. Let $\psi: \mathbb{D} \to \overline{\mathbb{C}}$, $\psi(0) = 0$, be a holomorphic ramified covering of local degree 3 over each point of $H$ and local degree 1 (unramified) over every point of $\overline{\mathbb{C}} \backslash H$. Such a ramified covering exists by the Uniformization Theorem for two-dimensional orbifolds [15, Ch. IX, Theorem 11]. Then $\psi$ has no asymptotic values and

$$\|\psi'(0)\| \geq 1 \tag{16}$$

by Schwarz's lemma, because $\psi$ is unramified over $\mathbb{D}$.

Now we set

$$g := h_2 \circ h_1 \circ \psi = h \circ \psi. \tag{17}$$

First we verify the statement about asymptotic values. Neither $h$ nor $\psi$ have them, so the composition $g$ does not.

The inequality (12) follows from (14), (16) and (17).

Now we verify

$$\mathfrak{B}(g) \leq \mathfrak{B}(p) + \varepsilon. \tag{18}$$

To prove (18) we assume that $B \subset \overline{\mathbb{C}}$ is a spherical disc of radius $R > \varepsilon$, where a branch of $g^{-1}$ is defined. By (17) there is a simply-connected region $D \subset \overline{\mathbb{C}}$ such that $h: D \to B$ is a homeomorphism, and a branch of $\psi^{-1}$ exists in $D$. It follows that $D \cap H = \emptyset$ so, by definition (15) of $H$

$$h_1(D) \cap E = \emptyset. \tag{19}$$

We consider the spherical disc $B_1 \subset B$ of spherical radius $R - \varepsilon$ and the same center as $B$. Let $D_1$ be the component of $h^{-1}(B_1)$ such that $D_1 \subset D$. Since $B$ is the open $\varepsilon$-neighborhood of $B_1$, and $h_2: h_1(D) \to B$ is an isometry, it follows that $h_1(D)$ is the open $\varepsilon$-neighborhood of $h_1(D_1)$. This, (19) and the definition of $E$ implies that $D_1 \subset \mathbb{D}$.

Since $D_1 \subset \mathbb{D}$, and $p$ is the restriction of $h$ on $\mathbb{D}$, $p$ maps $D_1$ onto $B_1$ homeomorphically, and so (18) follows. Together with (13) this gives (11).

It remains to check that all critical points of $g$ are multiple if all critical points of $f$ are. Using (17) we see that if $z_0$ is a critical point of $g$, then either $z_0$ is a critical point of $\psi$ or $\psi(z_0)$ is a critical point of $h$. In the first case $z_0$ is multiple, since all critical points of $\psi$ are multiple according to the definition of $\psi$. In the second case, when $z_0$ is not a critical point of $\psi$ and $\psi(z_0)$ is a critical point of $h$, we must have $\psi(z_0) \in \mathbb{D}$, since $\psi$ is ramified over all critical points of $h$ outside $\mathbb{D}$ by (15). But in $\mathbb{D}$ the maps $p$ and $h$ are the same, and $p$ has only multiple critical points. Hence $\psi(z_0)$ is a multiple critical point of $h$ which implies that $z_0$ is a multiple critical point of $g$. □

LEMMA 2.3 ([10, Prop. 8.4]). *Let $D$ be a Riemann surface, and $f\colon D \to \overline{\mathbb{C}}$ a nonconstant holomorphic map without asymptotic values such that $\mathfrak{B}(f) < \pi/2$.*

*Then there exists a set $T$ of compact topological triangles in $D$ with the following properties*:

(a) *For all $\Delta \in T$, the edges of $\Delta$ are analytic arcs, $f|\Delta$ is injective and conformal on $\mathring{\Delta}$. The set $f(\Delta)$ is a spherical triangle contained in a closed spherical disc of radius $\mathfrak{B}(f)$.*

(b) *If the intersection of two distinct triangles $\Delta_1, \Delta_2 \in T$ is nonempty, then $\Delta_1 \cap \Delta_2$ is a common edge of $\Delta_1$ and $\Delta_2$ or a set of common vertices.*

(c) *The set consisting of all the vertices of the triangles $\Delta \in T$ is equal to the set of critical points of $f$.*

(d) *$T$ is locally finite, i.e., for every $z \in D$ there exists a neighborhood $W$ of $z$ such that $W \cap \Delta \neq \emptyset$ for only finitely many $\Delta \in T$.*

(e) $\bigcup_{\Delta \in T} \Delta = D$.

A complete proof of this Lemma 2.3 is contained in [10, §8]. Here we just sketch the construction for the reader's convenience. Take a point $z \in D$ such that $f'(z) \neq 0$ and put $w = f(z)$. Then there exists a germ $\phi_z$ of $f^{-1}$ such that $\phi_z(w) = z$. Let $B \subset \overline{\mathbb{C}}$ be the largest open spherical disc centered at $w$ to which $\phi_z$ can be analytically continued. As there are no asymptotic values, the only obstacle to analytic continuation comes from the critical points of $f$. So there is at least one but at most finitely many singularities of $\phi_z$ on the boundary $\partial B$. Let $C(z)$ be the spherical convex hull of these singularities. This is a spherically convex polygon contained in $\overline{B}$. This polygon is nondegenerate (has nonempty interior) if and only if the number of singular points on $\partial B$ is at least three. Let $D(z) = \phi_z(C(z))$. We consider the set $Q$ of all points $z$ in $D$, for which the polygon $C(z)$ is nondegenerate. It can be shown that the union of the sets $D(z)$ over all $z \in Q$ is a locally finite covering of $D$, and $D(z_1) \cap D(z_2)$ for two different points $z_1$ and $z_2$ in $Q$ is either empty or a common side or a set of common vertices. The vertices are exactly the critical points of $f$. Finally, if we partition each $C(z)$, $z \in Q$, into spherical triangles by drawing appropriate diagonals, then the images of these triangles under the maps $\phi_z$, $z \in Q$, give the set $T$. □

*Reduction of Theorems* 1.2 *and* 1.3 *to Theorem* 1.4. We assume that $R = 1$ in Theorems 1.2 and 1.3. This does not restrict generality because we can replace $f(z)$ by $f(z/R)$. Moreover, the hypotheses on the function $f$ in both theorems and the estimates (1) and (3) possess an obvious invariance with respect to pre-composition of $f$ with an automorphism of the unit disc. So it is enough to prove (1) and (3) for $z = 0$.

By Lemma 2.2 it will suffice to consider only the case when $f$ has no asymptotic values. If $\mathfrak{B}(f) \geq \pi/2$ there is nothing to prove, so we assume that $\mathfrak{B}(f) < \pi/2$. Now we apply Lemma 2.3 to $f: \mathbb{D} \to \overline{\mathbb{C}}$. As it was explained in the Introduction, we equip $\mathbb{D}$ with the pullback of the spherical metric via $f$, which turns $\mathbb{D}$ into an open simply-connected singular surface, which we call $S_f$. The map $f: \mathbb{D} \to \overline{\mathbb{C}}$ now factors as in (5), where $\mathbb{D}$ is the unit disc with the usual hyperbolic metric, $f_1$ is a homeomorphism (coming from the identity map), and $f_2$ is a path isometry. The restriction of $f_2$ onto each triangle $\Delta \in T$ is an isometry, so we can call $\Delta$ a spherical triangle. The circumscribed radius of each triangle $\Delta \in T$ is at most $\mathfrak{B}(f)$, and the sum of angles at each vertex is at least $4\pi$. Under the hypotheses of Theorem 1.3 it is at least $6\pi$. Thus one of the conditions (i) or (ii) of Theorem 1.4 is satisfied.

As $f$ has no asymptotic values the singular surface $S_f$ is complete. Applying Theorem 1.4 to $S = S_f$, we find a $K$-quasiconformal and $L$-Lipschitz homeomorphism $\psi^{-1}: \mathbb{D} \to S_f$ with $K \geq 1$ and $L \geq 1$ depending only on $\varepsilon$.

Now we put $\phi = \psi \circ f_1: \mathbb{D} \to \mathbb{D}$. This map $\phi$ is a $K$-quasiconformal homeomorphism. Post-composing $\psi$ with a conformal automorphism of $\mathbb{D}$, we may assume $\phi(0) = 0$. Then Mori's theorem [4, IIIC] yields that $|\phi(z)| \leq 16|z|^{1/K}$ for $z \in \mathbb{D}$. This implies

$$\operatorname{dist}_h(0, \phi(z)) \leq 43|r|^{1/K} \quad \text{for} \quad z \in D(r),\ r \in (0, 32^{-K}],$$

where $\operatorname{dist}_h$ denotes the hyperbolic distance. Thus by the Lipschitz property of $\psi^{-1}$ we obtain that $f_1 = \psi^{-1} \circ \phi$ maps $D(r)$, $r \in (0, 32^{-K}]$, into a disc on $S_f$ centered at $f_1(0)$ of radius at most $43Lr^{1/K}$. Since $f_2$ is a path isometry, we conclude that $f = f_2 \circ f_1$ maps the disc $D(r_0)$ with $r_0 = (32L)^{-K}$ into a hemisphere centered at $f(0)$. Then Schwarz's lemma implies $\|f'(0)\| \leq (32L)^K$. □

It remains to prove Theorem 1.4. We begin with the study of a projection map $\Pi$ which associates a Euclidean triangle with each spherical triangle.

LEMMA 2.4. *Let $\Delta \subset \Sigma$ be a spherical triangle of spherical circumscribed radius $R < \pi/2$, $C$ its circumscribed circle, and $P \subset \mathbb{R}^3$ the plane containing $C$. Let $\Pi: \Delta \to P$ be the central projection from the origin. Then $\Pi$ is an $L$-bilipschitz map from $\Delta$ onto $\Pi(\Delta)$ with $L = \sec R$. Furthermore the ratios of the angles of $\Pi(\Delta)$ to the corresponding angles of $\Delta$ are between* $\cos R$ *and* $\sec R$.

Further properties of the map $\Pi$ are stated in Lemmas 3.1 and 3.5. Note that the triangle $\Pi(\Delta)$ is congruent to $F\Delta$, where $F = \mathrm{chd}$.

*Proof of Lemma* 2.4. Take any point $p \in \Delta$ and consider the tangent plane $P_1$ to $\Sigma$ passing through $p$. The angle $\gamma < \pi/2$ between $P_1$ and $P$ is at most $R$. Suppose that $p$ is different from the spherical center $p_0$ of $C$. Then there is a unique unit vector $u$ tangent to $\Sigma$ at $p$ pointing towards $p_0$. Let $v$ be a unit vector perpendicular to $u$ in the same tangent plane. Then the maximal length distortion of $\Pi$ at $p$ occurs in the direction of $u$ and is $\cos R \sec^2 \gamma$. The minimal length distortion occurs in the direction of $v$ and is $\cos R \sec \gamma$. These expressions for the maximal and minimal length distortion of $\Pi$ are also true at $p = p_0$ where $\gamma = 0$. Considering the extrema of the maximal and minimal distortion for $\gamma \in [0, R]$ we obtain $L = \sec R$.

To study the angle distortion we consider a vertex $O$ of $\Delta$. Let $v$ be a unit tangent vector to $P$ at $O$, and its preimage $u = (\Pi^{-1})'(v)$ in the tangent plane to $\Sigma$ at $O$. If $v$ makes an angle $\tau \in (0, \pi)$ with $C$, then the components of $v$ tangent and normal to $C$ have lengths $|\cos \tau|$ and $\sin \tau$, respectively. The component tangent to $C$ remains unchanged under $(\Pi^{-1})'$, while the normal one decreases by the factor of $\cos R$. Thus the angle between $u$ and $C$ is

$$\eta = \operatorname{arccot}(\sec R \cot \tau). \tag{20}$$

This distortion function has derivative

$$\frac{d\eta}{d\tau} = \frac{\cos R}{\cos^2 \tau + \cos^2 R \sin^2 \tau}, \tag{21}$$

which is increasing from $\cos R$ to $\sec R$ as $\tau$ runs from $0$ to $\pi/2$. □

LEMMA 2.5. *Suppose that $\phi: \Delta \to \mathbb{C}$ is an affine map of a Euclidean triangle, whose angles are $\alpha \le \beta \le \gamma$. Let $L$ be the Lipschitz constant of $\phi$, and let $l$ be the maximal Lipschitz constant of the three maps obtained by restricting $\phi$ to the sides of $\Delta$. Then $L \le l\pi^2\beta^{-2}$.*

*Proof.* Let $u_1$ and $u_2$ be unit vectors in the directions of those sides of $\Delta$ that form the angle $\beta$. Let $u$ be a unit vector for which

$$L = |\phi' u|. \tag{22}$$

If $u = c_1 u_1 + c_2 u_2$, then by taking scalar products we obtain

$$\begin{aligned}(u, u_1) &= c_1 + c_2(u_1, u_2),\\ (u, u_2) &= c_1(u_1, u_2) + c_2.\end{aligned}$$

Solving this system with respect to $c_1$ and $c_2$ by Cramer's rule, and using trivial estimates for products, we get

$$|c_j| \leq \frac{2}{1-(u_1,u_2)^2} = 2\csc^2(\beta).$$

By substituting this into (22), we conclude $L \leq (|c_1| + |c_2|)l \leq 4l\csc^2(\beta) \leq l\pi^2\beta^{-2}$. □

Lemma 2.6. *For $q \in (0,1)$ the metric in the disc $D(\sqrt{2})$ given by the length element*

$$\lambda(z)|dz| := \frac{2qR^q|z|^{q-1}|dz|}{R^{2q}-|z|^{2q}}, \quad \textit{where} \quad R := \sqrt{2}\left(\frac{1+q}{1-q}\right)^{1/(2q)}, \tag{23}$$

*has constant curvature $-1$ everywhere in $D(\sqrt{2})\backslash\{0\}$. Its density $\lambda$ is a decreasing function of $|z| \in (0,\sqrt{2})$ with infimum $\sqrt{(1-q^2)/2}$.*

This is proved by a direct computation. □

We will repeatedly use the following facts.

If $F$ is a concave nondecreasing function with $F(0) = 0$ and $F(x) > 0$ for $x > 0$, then the inequalities $a, b > 0$ and $c < a + b$ imply $F(c) < F(a) + F(b)$. Thus for every spherical or Euclidean triangle $\Delta$ the transformed Euclidean triangle $\tilde{\Delta} = F\Delta$ is defined.

This applies to both side distortion functions in (8) and (9).

If $\alpha \leq \beta \leq \gamma$ are the angles of $\Delta$, then the corresponding angles $\tilde{\alpha}, \tilde{\beta}, \tilde{\gamma}$ of $\tilde{\Delta}$ satisfy $\tilde{\alpha} \leq \tilde{\beta} \leq \tilde{\gamma}$. This follows from a well-known theorem of elementary geometry that larger angles be opposite larger sides.

*Derivation of Theorem* 1.4 *from Theorems* 1.5 *and* 1.6. Let us assume that condition (i) of Theorem 1.4 holds. (The proof under condition (ii) is similar.)

We will construct the required map in two steps. First we will map the given singular surface $S$ onto a singular surface $\tilde{S}$ which has the flat Euclidean metric everywhere except at a set consisting of isolated singularities, where we have some definite positive total angle excess.

To each triangle $\Delta \in T$ we assign a Euclidean triangle $\tilde{\Delta} = F_k \Delta$. The subadditivity of $F_k$ ensures that $\tilde{\Delta}$ is well-defined. If $\Delta_1$ and $\Delta_2$ in $T$ have a common side or common vertices, then we identify the corresponding sides or vertices of $\tilde{\Delta}_1$ and $\tilde{\Delta}_2$. For the identification of common sides we use arclength as the identifying function. By gluing the triangles $\tilde{\Delta}$ together in this way, we obtain a new singular surface $\tilde{S}$. It is Euclidean everywhere except at the vertices of the triangles $\tilde{\Delta}$. The total angle at each vertex is at least $2\pi(1+2\delta)$ by Theorem 1.5.

Now we construct a bilipschitz map $\psi_1: S \to \tilde{S}$. We will define it on each triangle in $T$ in such a way that the definitions match on the common sides and vertices of the triangles. For a given triangle $\Delta \in T$ we put $\psi_1|_\Delta = \phi \circ \Pi$, where $\Pi$ is the central projection map from Lemma 2.4, and $\phi$ is the unique affine map of the Euclidean triangles $\Pi(\Delta) \to \tilde{\Delta}$.

According to the identifications used to define $\tilde{S}$ and since our maps between triangles map vertices to vertices, it is clear that the definition of $\psi_1$ matches for common vertices of triangles. Let $s$ be a common side of two triangles $\Delta_1$ and $\Delta_2$ in $T$. Then we place $\Delta_1$ and $\Delta_2$ on the sphere $\Sigma$ in such a way that they have this common side and consider the planes $P_1$ and $P_2$ in $\mathbb{R}^3$ passing through the vertices of $\Delta_1$ and $\Delta_2$, respectively. Then the central projections $\Pi_1$ and $\Pi_2$ from $\Sigma$ to $P_1$ and $P_2$, respectively, match on $s$, which is mapped by both projections onto the chord of $\Sigma$ connecting the endpoints of $s$ in $\mathbb{R}^3$. That the affine maps $\phi_1: \Pi_1(s) \to \tilde{S}$ and $\phi_2: \Pi_2(s) \to \tilde{S}$ match is evident.

Our surfaces $S$ and $\tilde{S}$ carry intrinsic metrics. Therefore, in order to show that $\psi_1$ is $L$-bilipschitz it suffices to show that the restriction of $\psi_1$ to an arbitrary triangle $\Delta$ in $T$ is $L$-bilipschitz. Since the circumscribed radius of $\Delta$ is bounded away from $\pi/2$, Lemma 2.4 gives a bound for the bilipschitz constant for the projection part $\Pi$ of $\psi_1$ independent of $\varepsilon$. (In case (ii) of Theorem 1.4 it will depend on $\varepsilon$.) To get an estimate for the affine factor $\phi$ we first consider length distortion on the sides of $\Pi(\Delta)$. According to the definition of $F_k$ in (8), a side of $\Pi(\Delta)$ of length $a$ is mapped onto a side of $\tilde{\Delta}$ with length $\min\{ka, \sqrt{a}\}$. Since $0 < a < 2$ and $k \geq 1$, we obtain $1/\sqrt{2} \leq \min\{ka, \sqrt{a}\}/a \leq k$. So the length distortion of $\phi$ on the sides is at most $\max\{k, \sqrt{2}\}$. Moreover, we note that the diameter of each triangle $\tilde{\Delta}$ is less than $\sqrt{2}$. To estimate the distortion in the interior of $\Pi(\Delta)$ we consider two cases.

If the diameter $d$ of $\Pi(\Delta)$ is at most $1/k^2$, then the triangle $\tilde{\Delta}$ is obtained from $\Pi(\Delta)$ by scaling its sides by the factor $k$. In particular, these triangles are similar, so the affine map $\phi$ is a similarity, and the distortion in the interior of $\Pi(\Delta)$ is equal to the distortion on the sides.

To deal with the the case $d > 1/k^2$ first we note that the circumscribed radius $r$ of $\Pi(\Delta)$ is at most 1. If we denote by $\alpha' \leq \beta'$ the two smaller angles of $\Pi(\Delta)$, then $\alpha' + \beta' \geq \arcsin(d/(2r)) \geq 1/(2k^2)$, and so $\beta' \geq 1/(4k^2)$. Then Lemma 2.5 gives an estimate of the Lipschitz constant of $\phi$. To estimate the Lipschitz constant of $\phi^{-1}$ we notice that the intermediate angle $\tilde{\beta}$ of $\tilde{\Delta}$ satisfies $\tilde{\beta} > \beta/2 \geq (\beta' \cos R)/2$ in view of Theorem 1.5 and Lemma 2.4. This gives $\tilde{\beta} \geq \cos b_0/(8k^2)$. So Lemma 2.5 gives a bound for the Lipschitz constant of $\phi^{-1}\colon \tilde{\Delta} \to \Pi(\Delta)$ as well.

In any case, we see that $\psi_1$ restricted to any triangle in $T$ is $L$-bilipschitz with bilipschitz constant only depending on $k$, and hence only on $\varepsilon$. As we stated above this implies that $\psi_1\colon S \to \tilde{S}$ is $L$-bilipschitz. Then $\psi_1$ is also $K$-quasiconformal with $K = L^2$.

Thus we have proved that $\psi_1\colon S \to \tilde{S}$ is $L$-bilipschitz and $K$-quasiconformal with $L$ and $K$ depending only on $\varepsilon$.

Now we proceed to the second step of our construction, and find a conformal map $\psi_2\colon \tilde{S} \to \mathbb{D}$. Since $S$ is an open and simply-connected surface, $\tilde{S}$ has the same properties. By the Uniformization Theorem there exists a conformal map $g\colon D(R) \to \tilde{S}$, where $0 < R \leq \infty$. We will estimate $\|g'\|$. If the total angle at a vertex $v \in \tilde{S}$ is

$$2\pi\alpha \geq 2\pi(1 + 2\delta) \tag{24}$$

and $z_0 = g^{-1}(v)$, then we have

$$|g'(z)| \sim \mathrm{const}|z - z_0|^{\alpha - 1}, \quad z \to z_0, \tag{25}$$

and

$$\rho_v(g(z)) \sim \mathrm{const}|z - z_0|^{\alpha}, \quad z \to z_0, \tag{26}$$

where $\rho_v(w)$ stands for the distance from a point $w \in \tilde{S}$ to the vertex $v \in \tilde{S}$. Now we put a new conformal metric on $\tilde{S}$. Denote by $V$ the set of all vertices of our covering $\{\tilde{\Delta}\}$ of $\tilde{S}$. For a point $w \in \tilde{S}$ we put $\rho(w) := \inf\{\rho_v(w) : v \in V\}$, which is the distance from $w$ to the set $V$. The infimum is actually attained, because our singular surface $\tilde{S}$ is complete and thus there are only finitely many vertices within a given distance from any point $w \in \tilde{S}$. As we noticed above the diameter of each triangle $\tilde{\Delta} \subset \tilde{S}$ is less than $\sqrt{2}$. Hence $\rho(w) < \sqrt{2}$ for $w \in \tilde{S}$. Let $\lambda$ be the density in (23) with $q := (1+\delta)^{-1}$. Then (24) implies that for every vertex with total angle $2\pi\alpha$ we have

$$\alpha q > 1. \tag{27}$$

Following Ahlfors we define a conformal length element $\Lambda(w)|dw|$ with the density

$$\Lambda(w) := \lambda(\rho(w)), \quad w \in \tilde{S} \backslash V.$$

Since $\rho < \sqrt{2}$ this is well-defined. For each point $p \in \tilde{S}\backslash V$ we can choose a vertex $v(p) \in V$ closest to $p$. Then in a neighborhood of $p$ the metric $\lambda(\rho_{v(p)}(w))|dw|$ is a supporting metric of curvature $-1$ in the sense of [5, §1-5] for $\Lambda(w)|dw|$. In view of (25), (26), (23) and (27) the density of the pullback of the metric $\Lambda(w)|dw|$ via the map $g$ has the following asymptotics near the preimage $z_0$ of a vertex $v$

$$\Lambda(g(z))|g'(z)| = \lambda(\rho_v(g(z)))|g'(z)| = O(|z-z_0|^{\alpha q-1}) = o(1), \quad z \to z_0.$$

The Ahlfors-Schwarz lemma and Lemma 2.6 now imply that for arbitrary $0 < r < R$

$$|g'(z)| \le (\inf \Lambda)^{-1} \frac{2r}{r^2 - |z|^2} \le \frac{2\sqrt{2}r}{\sqrt{1-q^2}(r^2-|z|^2)}, \quad z \in D(r). \tag{28}$$

This inequality shows that $R < \infty$; thus we can assume without loss of generality that $R = 1$. Then (28) is true for $r = 1$, and this implies $\|g'\| \le \sqrt{2}(1-q^2)^{-1/2}$. In other words, if we put $\psi_2 := g^{-1}$, then $\psi_2^{-1}$ is a Lipschitz map with Lipschitz constant $\sqrt{2}(1-q^2)^{-1/2}$ which depends only on $\delta$ and hence only on $\varepsilon$.

Composing our maps we obtain $\psi = \psi_2 \circ \psi_1 : S \to \mathbb{D}$. Then $\psi$ is $K$-quasiconformal and $\psi^{-1}$ is $L$-Lipschitz with $K$ and $L$ depending only on $\varepsilon$. This proves Theorem 1.4. □

## 3. Outline of the proof of Theorem 1.5. The generic case

We use the notation $F\Delta$ and $D(F, \Delta)$ defined in (7). For the proof of Theorem 1.5 we need first of all an increasing subadditive function $F$ with $F(0) = 0$ such that

$$D(F, \Delta) > 1/2 \tag{29}$$

for every spherical triangle $\Delta$ of circumscribed radius less than $b_0$. An analytic function $F$ with these properties is $F_\infty := \sqrt{\text{chd}}$. That $F_\infty$ indeed satisfies (29) is the core of our argument (Lemmas 3.1–3.4). We split $F_\infty$ into the composition of chd and $\sqrt{\ }$ and introduce the intermediate Euclidean triangle $\Delta' = \text{chd}\,\Delta$, which is obtained by replacing the sides of $\Delta$ by the corresponding chords. Then $\tilde{\Delta} := F_\infty \Delta = \sqrt{\ }\Delta'$. Replacing the sides by their chords may decrease the angles by a factor of $1/3$, and taking square roots of the sides of a Euclidean triangle may decrease the angles by a factor of $1/2$. Nevertheless, the two parts of our map somehow compensate each other, and we get (29) for $F = F_\infty$.

In fact the function $F_\infty$ is not good enough for our purposes for two reasons. First, $F'_\infty(0) = \infty$, so the map it induces cannot be bilipschitz.

Second, inequality (29) with $F = F_\infty$ is best possible, while we need a lower estimate of the form $1/2+\delta$. The distortion $D(F_\infty, \Delta)$ becomes close to $1/2$ in two situations. First, when $\Delta$ is close to the extremal equilateral triangle. This is avoided by our assumption that the circumscribed radius of a triangle is bounded away from $b_0$. Second, the distortion $D(F_\infty, \Delta)$ may tend to $1/2$ as $\Delta$ degenerates.

A natural remedy for both drawbacks of $F_\infty$ is to introduce a 'cutoff,' that is, to use $F_k = \min\{k \operatorname{chd}, \sqrt{\operatorname{chd}}\}$ with sufficiently large $k$. Unfortunately, for $F = F_k$ it is harder to prove (29) directly.

We split all spherical triangles into several categories, according to the degree of their degeneracy. Triangles with only one small angle are dealt with in Section 4. This is a relatively simple case. The worst problems occur when a triangle has two small angles. But given that the circumscribed radius is bounded by $b_0 < \pi/2$, such triangles have to be small (Lemma 3.5). So their distortion can be studied using Euclidean trigonometry. This we do in Section 5. Finally, in Section 6 we derive Theorem 1.5 from all these special cases.

The proof of Theorem 1.6 consists of similar steps, but each step is much easier. For instance, the counterparts of the Lemmas 3.2–3.4 and 5.3 are completely trivial.

Now we proceed to the formal proofs.

In the rest of this section we consider spherical triangles $\Delta \subset \Sigma$ with one marked vertex $O$. The sides of $\Delta$ are $a, b, c$ and the angles opposite to these sides are $\alpha, \beta, \gamma$. We assume that the angle at $O$ is $\gamma$.

We first study the distortion of the angle $\gamma$ under the transformation $\tilde{\Delta} = F_\infty \Delta$, defined in the Introduction, where

$$F_\infty(t) = \sqrt{\operatorname{chd} t}. \tag{30}$$

For this purpose we do the following construction. Let $C \subset \Sigma$ be the circumscribed circle of $\Delta$ and $r$ its Euclidean radius. Let $R < \pi/2$ be the spherical circumscribed radius of $\Delta$, that is, the radius of the spherical disc in $\Sigma$ that contains $\Delta$ and is bounded by $C$. Then

$$r = \sin R. \tag{31}$$

As in Lemma 2.4, $P$ is the plane containing the vertices of $\Delta$ (and thus containing $C$), and $\Pi$ is the central projection of the circumscribed disc to $P$ from the origin. Then $\Delta' = \Pi(\Delta) \subset P$ is a Euclidean triangle with sides $\operatorname{chd} a, \operatorname{chd} b, \operatorname{chd} c$, and $\tilde{\Delta} = F_\infty \Delta = G_\infty \Delta'$, where

$$G_\infty(x) = \sqrt{x}. \tag{32}$$

We use the same letters to denote the corresponding sides and angles of the triangles $\Delta, \Delta'$ and $\tilde{\Delta}$ with primes for $\Delta'$ and tildes for $\tilde{\Delta}$.

Let $2\phi = \gamma'$, $\phi \in (0, \pi/2)$, and $t \in (\phi, \pi - \phi)$ be the angle between the bisector of $\gamma'$ and the circle $C$ at the point $O$. Then the three parameters $(R, \phi, t)$ describe all possible spherical triangles with a chosen vertex. See Figure 1 below.

First we study the angle distortion by $\Pi$.

Lemma 3.1. *Among all spherical triangles $\Delta$ with fixed circumscribed radius $R < \pi/2$ and fixed angle $\gamma' = 2\phi$ the largest angle $\gamma$ occurs when $\Delta$ is isosceles with two equal sides meeting at $\gamma$.*

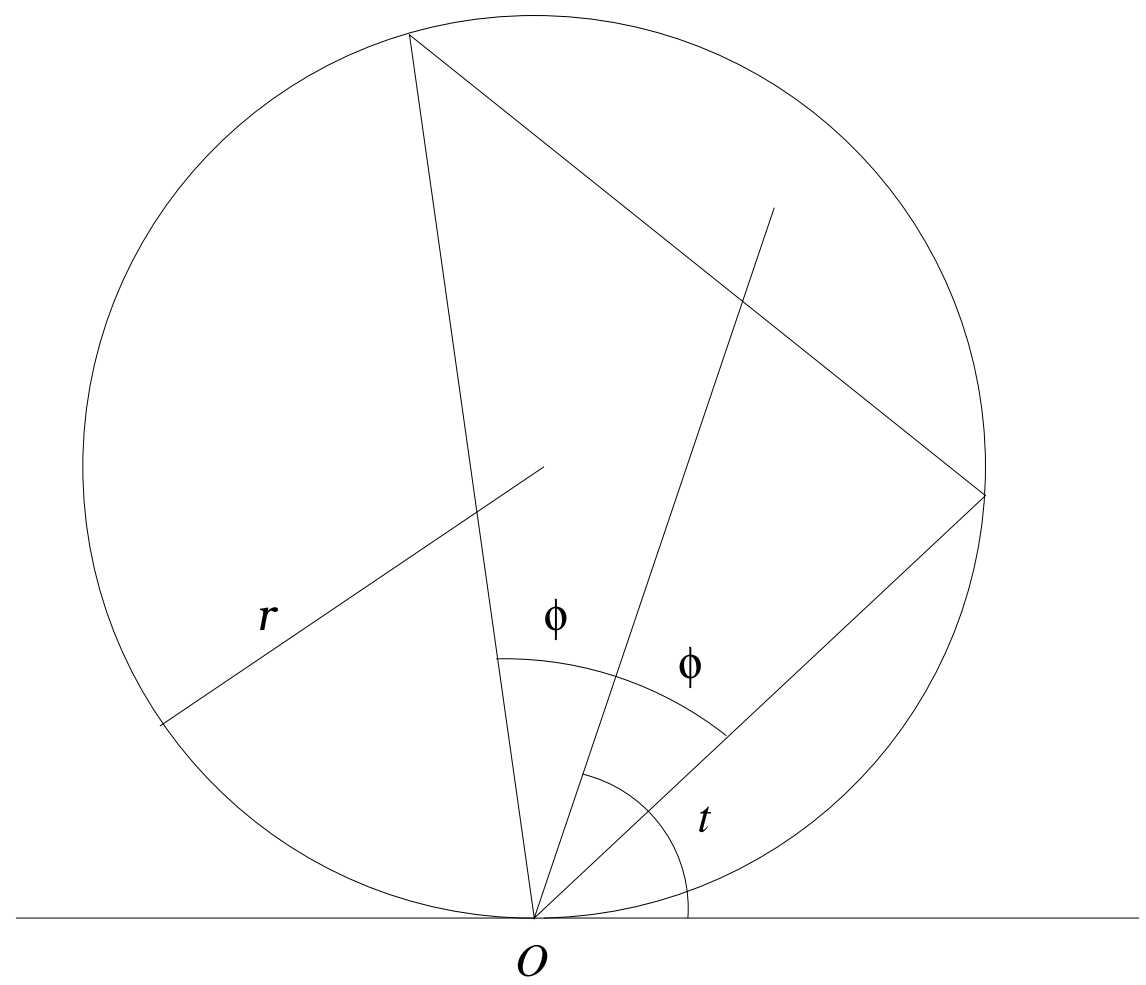

Figure 1. Plane $P$

*Proof.* From (20) and (21) we obtain

$$(33) \quad \gamma = \gamma(t) = \operatorname{arccot}(\sec R \cot(t+\phi)) - \operatorname{arccot}(\sec R \cot(t-\phi)) = \int_{t-\phi}^{t+\phi} \frac{\cos R \, d\tau}{\cos^2 \tau + \cos^2 R \sin^2 \tau} =: \int_{t-\phi}^{t+\phi} \mu(\tau)\, d\tau.$$

Then $\mu$ is increasing on $[0, \pi/2]$ and satisfies $\mu(\tau) = \mu(\pi/2 - \tau)$. So for fixed $\phi$ the maximum of the integral occurs when $t = \pi/2$, which proves the lemma. □

Now we consider the distortion of the angle $\gamma'$ of the Euclidean triangle $\Delta' = \Pi(\Delta)$ under the transformation $\Delta' \mapsto G_\infty \Delta' = \tilde{\Delta}$ with $G_\infty$ as in (32).

Lemma 3.2. *Among all Euclidean triangles $\Delta'$ with fixed circumscribed radius $r$ and fixed angle $\gamma' = 2\phi$ the smallest angle $\tilde{\gamma}$ in $\tilde{\Delta} = G_\infty \Delta'$ occurs when $\Delta'$ is isosceles with two equal sides meeting at $\gamma'$.*

*Proof.* Since the triangles we consider are Euclidean and the distortion function $G_\infty$ is homogeneous, we may assume without loss of generality that $r = 1$. Then the sides of $\Delta'$ are (see Figure 1)

$$a' = 2\sin(t-\phi), \quad b' = 2\sin(t+\phi), \quad c' = 2\sin 2\phi,$$

and the sides of $\tilde{\Delta}$ are $\sqrt{a'}$, $\sqrt{b'}$, $\sqrt{c'}$. By the Law of Cosines

$$\begin{aligned}\cos\tilde{\gamma}(t) &= \frac{\sin(t-\phi)+\sin(t+\phi)-\sin 2\phi}{2\sqrt{\sin(t-\phi)\sin(t+\phi)}} \\ &= \frac{2\sin t\cos\phi - 2\sin\phi\cos\phi}{2\sqrt{\sin^2 t - \sin^2\phi}} \\ &= \cos\phi\sqrt{1-\frac{2\sin\phi}{\sin t+\sin\phi}}.\end{aligned}$$

The last expression has a maximum for $\phi \le t \le \pi/2$ at $t = \pi/2$. □

As a corollary from Lemmas 3.1 and 3.2 we obtain the following.

Lemma 3.3. *Among all spherical triangles $\Delta$ with fixed circumscribed radius $R < \pi/2$ and fixed angle $\gamma'$ the ratio $\tilde{\gamma}/\gamma$ is minimal for the isosceles triangle $\Delta$ with two equal sides meeting at $\gamma$. Here $\tilde{\gamma}$ is the angle of $\tilde{\Delta} = F_\infty\Delta$, corresponding to $\gamma$.*

*Proof.* It is clear that $\Delta' = \Pi(\Delta)$ is isosceles with two equal sides meeting at $\gamma'$ if and only if $\Delta$ is isosceles with two equal sides meeting at $\gamma$. Using Lemmas 3.1 and 3.2 and the notation introduced in their proofs, we have $\gamma'(t) \equiv 2\phi = \text{const}$ for $\phi \le t \le \pi/2$; thus

$$\frac{\tilde{\gamma}(t)}{\gamma(t)} \ge \frac{\tilde{\gamma}(\pi/2)}{\gamma(\pi/2)}.$$ □

Now we state and prove the main result of this section.

Lemma 3.4. *If $F_\infty$ is the function of* (30), *then the distortion satisfies $D(F_\infty,\Delta) > 1/2$ for all spherical triangles $\Delta$ of circumscribed radius $R < b_0$.*

*Proof.* In view of Lemma 3.3 it is enough to show that for every isosceles spherical triangle $\Delta$ of circumscribed radius $R < b_0$ which has two equal sides meeting at $\gamma$, we have $\tilde{\gamma} > \gamma/2$. In other words, we fix $R < b_0$ and $t = \pi/2$ and study the ratio $\tilde{\gamma}/\gamma$ as a function of $\phi$.

We have

$$\operatorname{chd} a = \operatorname{chd} b = r \operatorname{chd}(\pi - 2\phi) = 2r\cos\phi \tag{34}$$

and $\operatorname{chd} c = r \operatorname{chd} 4\phi = 4r\cos\phi\sin\phi$, where $r$ is the same as in (31). Eliminating $\phi$ from these relations and using (31) we obtain

$$\sin\frac{c}{2} = 2\sin\frac{a}{2}\sqrt{1-\frac{\sin^2(a/2)}{\sin^2 R}}. \tag{35}$$

Now the lengths of the sides of $\tilde{\Delta}$ are the square roots of the lengths of the sides $\operatorname{chd} a, \operatorname{chd} a, \operatorname{chd} c$ of $\Delta'$, so by the Law of Cosines and (35)

$$\cos\tilde{\gamma} = \frac{2\operatorname{chd} a - \operatorname{chd} c}{2\operatorname{chd} a} = 1 - \frac{\sin(c/2)}{2\sin(a/2)} = 1 - \sqrt{1-\frac{\sin^2(a/2)}{\sin^2 R}}. \tag{36}$$

For the angle $\gamma$ we use (33) to obtain

$$\gamma = \pi - 2\arctan\left[\cos R\tan(\pi/2-\phi)\right] =: \pi - 2\arctan U.$$

So, using (34) and (31)

$$\begin{aligned} \cos\frac{\gamma}{2} &= \sin(\arctan U) = \frac{U}{\sqrt{1+U^2}} \\ &= \frac{\cos R\cot\phi}{\sqrt{1+\cos^2 R\cot^2\phi}} \\ &= \frac{\cos R\cos\phi}{\sqrt{\sin^2\phi+\cos^2 R\cos^2\phi}} \\ &= \frac{\cot R\operatorname{chd} a}{2\sqrt{1-\sin^2 R\cos^2\phi}} \\ &= \frac{\cot R\sin(a/2)}{\cos(a/2)} = \frac{\tan(a/2)}{\tan R}. \end{aligned} \tag{37}$$

Our goal is to prove that $\tilde{\gamma} > \gamma/2$, which is equivalent to

$$\cos\tilde{\gamma} < \cos(\gamma/2). \tag{38}$$

In view of (37) and (36) this is equivalent to

$$1-\sqrt{1-\frac{\sin^2(a/2)}{\sin^2 R}} < \frac{\tan(a/2)}{\tan R}.$$

Introducing temporary notation $T = \tan R \in (0,\infty)$ and $t = \tan(a/2) \in (0,T]$ we rewrite the previous inequality as

$$1-\frac{t}{T} < \sqrt{1-\frac{1+T^2}{T^2}\frac{t^2}{1+t^2}} = \frac{\sqrt{T^2-t^2}}{T\sqrt{1+t^2}},$$

or

$$\sqrt{(T-t)(1+t^2)} < \sqrt{T+t},$$

which after simplification becomes

$$t^2 - Tt + 2 > 0.$$

This holds for every $t \in (0, T]$ if $T^2 - 8 < 0$, which is $\tan R < 2\sqrt{2}$, and this is the same as $R < b_0$. This proves (38) and Lemma 3.4. □

This is the only place in our proof where the numerical value of $b_0$ is essential.

The following lemma is used in our study of degenerating triangles, which will occupy the rest of the proof.

LEMMA 3.5. *Let $\Delta$ be a spherical triangle of circumscribed radius $R < \pi/2$, and spherical diameter $d$, let $\alpha \leq \beta \leq \gamma$ be its angles, and $\alpha' \leq \beta' \leq \gamma'$ be the corresponding angles of $\Delta' = \Pi(\Delta)$. Then*

$$(1/2) \operatorname{chd} d \cot R \leq \alpha + \beta \leq 4 \operatorname{chd} d \csc(2R). \tag{39}$$

*Moreover, if* $\operatorname{chd} d \leq 2^{-1/2} \sin(R/2)$, *then*

$$\alpha \leq \alpha' \quad \textit{and} \quad \beta \leq \beta'. \tag{40}$$

Notice that if $R$ is bounded away from $\pi/2$, then Lemma 3.5 implies $\alpha + \beta \asymp d/R$. Furthermore, $\alpha \leq \alpha'$, $\beta \leq \beta'$ if $d/R$ is small enough.

*Proof.* The Law of Sines shows that if $d' = \operatorname{chd} d$ is the diameter of $\Delta'$ and $r = \sin R$ is its circumscribed radius, then

$$\sin(\alpha' + \beta') = d'/(2r).$$

Since $\alpha'$ and $\beta'$ are the two smaller angles of $\Delta'$, we have $\alpha' + \beta' \leq 2\pi/3$. Hence

$$\alpha' + \beta' \geq \arcsin\left(\frac{d'}{2r}\right),$$

with equality if $d'/(2r) < \sin(\pi/3) = \sqrt{3}/2$. Using $x \leq \arcsin x \leq \pi x/2$ for $0 \leq x \leq 1$ we obtain for $d'/(2r) < \sqrt{3}/2$

$$\frac{d'}{2r} \leq \alpha' + \beta' \leq \frac{2d'}{r}. \tag{41}$$

Here the left inequality is always true. But if $d'/(2r) \geq \sqrt{3}/2$, then the right inequality also holds since $\alpha' + \beta' \leq 2\pi/3$. Lemma 2.4 gives $\alpha'/\alpha \in [\cos R, \sec R]$ and similar bounds for $\tilde{\beta}/\beta$. From this and (41) the inequality (39) follows.

To prove (40) we recall the angle distortion function $\eta = \operatorname{arccot}(\sec R \cot \tau)$ from (20), whose derivative, given by (21), is increasing on $(0, \pi/2)$. The equation $d\eta/d\tau = 1$ has a unique solution on $[0, \pi/2]$, namely

$$\tau_0 = \arcsin(2^{-1/2} \sec(R/2)) \geq 2^{-1/2} \sec(R/2).$$

This means that for every sector of opening $\tau$ in the plane $P$ both of whose sides are within an angle of $\tau_0$ from the tangent to the circumscribed circle, the spherical angle $\eta$ corresponding to $\tau$ under the map $\Pi$ is smaller than $\tau$.

Suppose the largest angle $\gamma'$ in a Euclidean triangle is greater than $\pi/2$. Consider a vertex $O \in C$ of a smaller angle. Then this smaller angle is contained in a sector with vertex $O$ and opening $\pi - \gamma'$, which is determined by the tangent to $C$ at $O$ and the largest side of the triangle. Using this we see that (40) is satisfied if $\pi - \gamma' = \alpha' + \beta' \leq \tau_0 < \pi/2$. By (41) this holds if $\operatorname{chd} d \leq 2^{-1/2} \sin(R/2)$. □

## 4. The case of one small angle

For a fixed circumscribed radius $R < \pi/2$ there are two ways a triangle can become degenerate: it may have one small angle or two small angles. In this section we consider the first of these two cases.

LEMMA 4.1. *Let $(H_k)$ be a sequence of concave increasing functions on $[0, \pi)$ with the following properties*:

1. $H_k(0) = 0$ *and* $H_k(x) > 0$ *for* $x > 0$, $k \in \mathbb{N}$.
2. $H_k(x)/x \leq M(x) < \infty$ *for* $x > 0$, $k \in \mathbb{N}$.
3. *If* $(x_k)$ *is a sequence in* $(0, \pi)$ *with* $\lim x_k = 0$, *then* $H_k(x_k)/x_k \to \infty$.

*We consider a sequence $(\Delta_k)$ of spherical triangles with circumscribed radii bounded by $R < \pi/2$. Let the angles of $\Delta_k$ be $\alpha_k \leq \beta_k \leq \gamma_k$, and assume that $\alpha_k \to 0$, $\beta_k \to \beta_0 \in (0, \pi)$, $\gamma_k \to \gamma_0 \in (0, \pi)$.*

*If $\tilde{\alpha}_k \leq \tilde{\beta}_k \leq \tilde{\gamma}_k$ are the angles of the transformed triangle $\tilde{\Delta}_k = H_k \Delta_k$, then*:

(i) $\underline{\lim}\, \tilde{\alpha}_k/\alpha_k > 1/2$, *and*

(ii) $\underline{\lim}\, \tilde{\beta}_k/\beta_k > 1/2$.

*Proof.* We will omit the subscript $k$ in the proof, and tacitly assume that all limits are understood for $k \to \infty$. As usual we will denote the sides of the triangle $\Delta$ by $a \leq b \leq c$, and the corresponding sides of the triangle $\tilde{\Delta}$ by $\tilde{a} \leq \tilde{b} \leq \tilde{c}$.

By the Spherical Law of Sines [12, p. 150], we have

$$\frac{\sin a}{\sin b} = \frac{\sin \alpha}{\sin \beta}. \tag{42}$$

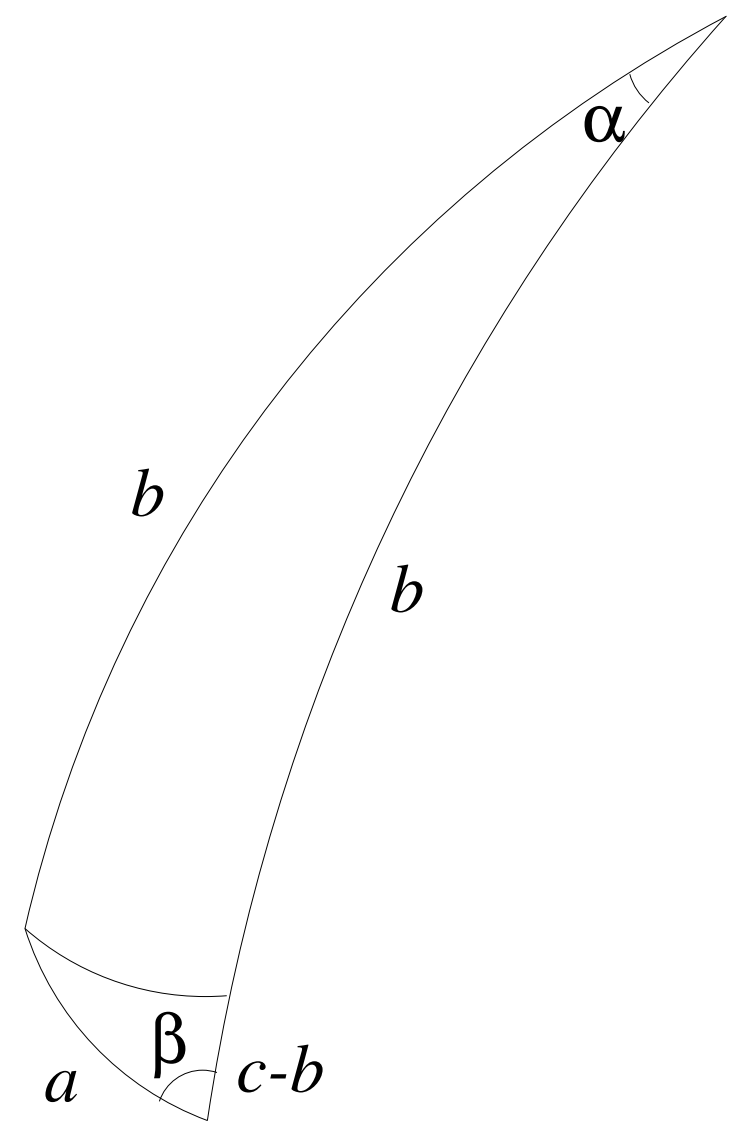

Figure 2. Triangle $\Delta$

By our assumption on the circumscribed radius of the triangles, the lengths of all sides are bounded away from $\pi$. Moreover, $\sin\beta \to \sin\beta_0 \neq 0$. Hence $\alpha \to 0$ implies $a \to 0$. Note that by the same reasoning we get $a/b \to 0$ which will be needed in the proof of Lemma 4.3.

Since $\Delta$ is contained in a sector of opening angle $\alpha \to 0$, we have

$$\operatorname{area}\Delta = \alpha + \beta + \gamma - \pi \to 0. \tag{43}$$

Hence $\beta + \gamma \to \pi = \beta_0 + \gamma_0$, which implies

$$\beta_0 \leq \pi/2. \tag{44}$$

One of Napier's Analogies [12, p. 160] shows

$$\frac{\tan\Big((c-b)/2\Big)}{\tan(a/2)} = \frac{\sin\Big((\gamma-\beta)/2\Big)}{\sin\Big((\gamma+\beta)/2\Big)} \to \cos\beta_0. \tag{45}$$

Since $a \to 0$ and all sides are bounded away from $\pi$, this gives $c - b \to 0$. But then we conclude from (45) that

$$\frac{c-b}{a} \to \cos\beta_0 \geq 0. \tag{46}$$

Concavity implies

$$\frac{H(c)-H(b)}{c-b} \leq \frac{H(b)}{b} \leq \frac{H(a)}{a}. \tag{47}$$

Using this, Mollweide's[3] formula (which is the Euclidean limit case of (45) [12, p. 61]), and (46) we get

$$\overline{\lim}\,\frac{\sin\Big((\tilde\gamma-\tilde\beta)/2\Big)}{\sin\Big((\tilde\gamma+\tilde\beta)/2\Big)} = \overline{\lim}\,\frac{\tilde c-\tilde b}{\tilde a}=\overline{\lim}\,\frac{H(c)-H(b)}{H(a)} \leq \lim\frac{c-b}{a}=\cos\beta_0. \tag{48}$$

On the other hand, since $\tilde\alpha+\tilde\beta+\tilde\gamma=\pi$ and $\tilde\alpha\leq\tilde\beta$ we have

$$\frac{\sin\Big((\tilde\gamma-\tilde\beta)/2\Big)}{\sin\Big((\tilde\gamma+\tilde\beta)/2\Big)}=\cos\tilde\beta-\tan(\tilde\alpha/2)\sin\tilde\beta\geq 2\cos\tilde\beta-1. \tag{49}$$

Since $\beta_0\leq\pi/2$ inequalities (48) and (49) imply

$$\overline{\lim}\,\cos\tilde\beta\leq(1+\cos\beta_0)/2=\cos^2(\beta_0/2)<\cos(\beta_0/2).$$

Here (and below in (50)) it is crucial that $\beta_0>0$. We get $\underline{\lim}\,\tilde\beta>\beta_0/2$ which implies (ii).

Since $\tilde\beta$ is the second largest angle in the Euclidean triangle $\tilde\Delta$, we have $\tilde\beta\leq\pi/2$. Therefore, from (44) and (ii) we get

$$\underline{\lim}\,\frac{\sin\tilde\beta}{\sin\beta}\geq\frac{\sin(\beta_0/2)}{\sin\beta_0}=(1/2)\sec(\beta_0/2)>1/2. \tag{50}$$

From the the Spherical and Euclidean Laws of Sines and the assumption that $\alpha\to 0$ we obtain

$$\underline{\lim}\,\tilde\alpha/\alpha \geq \underline{\lim}\,\frac{\sin\tilde\alpha}{\sin\alpha}=\underline{\lim}\,\frac{\tilde a\sin b\sin\tilde\beta}{\tilde b\sin a\sin\beta} \geq \underline{\lim}\,\frac{\sin\tilde\beta}{\sin\beta}\cdot\underline{\lim}\,\frac{\tilde a\sin b}{\tilde b\sin a}. \tag{51}$$

We consider

$$l:=\underline{\lim}\,\frac{\tilde a\sin b}{\tilde b\sin a}=\underline{\lim}\,\frac{H(a)}{a}\frac{b}{H(b)}\frac{\sin b}{b}. \tag{52}$$

Given any subsequence of our triangles, one can select from it another subsequence, such that either $\underline{\lim}\,b>0$ or $\lim b=0$. In the first case, the estimate $\sin x/x\asymp 1$ for $x>0$ bounded away from $\pi$, and properties 2 and 3 of the functions $H_k$ imply that $l=\infty$ in (52) for the selected subsubsequence. In the second case, when $\lim b=0$, we see from (47) that $l\geq 1$ in (52). This shows that for every subsequence $l\geq 1$. Now (i) follows from (51) and (50). □

We need a similar lemma for Euclidean triangles.

[3] Actually first published by Newton [21].

LEMMA 4.2. *Let $H$ be a concave increasing function on $[0,\infty)$ with $H(0)=0$ and $H(x)>0$ for $x>0$. We consider a sequence of Euclidean triangles $(\Delta_k)$ with angles $\alpha_k \le \beta_k \le \gamma_k$ and assume that $\alpha_k \to 0$, $\beta_k \to \beta_0 \in (0,\pi)$ and $\gamma_k \to \gamma_0 \in (0,\pi)$.*

*If $\tilde{\alpha}_k \le \tilde{\beta}_k \le \tilde{\gamma}_k$ are the angles of the transformed triangle $\tilde{\Delta}_k = H\Delta_k$, then*:

(i) $\underline{\lim}\, \tilde{\alpha}_k/\alpha_k > 1/2$, *and*

(ii) $\underline{\lim}\, \tilde{\beta}_k/\beta_k > 1/2$.

*Proof.* We follow the lines of the proof of Lemma 4.1 with the following modifications. We use the Euclidean Law of Sines instead of (42) to conclude that $a/b \to 0$. Formula (43) is replaced by the identity $\alpha+\beta+\gamma=\pi$ which implies $\beta+\gamma\to\pi$. Again (44) follows. Instead of Napier's analogy we use Mollweide's formula to obtain (46). Again (47) is true. Formulas (48), (49), and the proof of (ii) are exactly the same. The inequality (i) follows by a similar (and easier) computation as in (51) using the Euclidean Law of Sines, (50), and (47), namely,

$$\begin{aligned}\underline{\lim}\, \frac{\tilde{\alpha}}{\alpha} \quad &\ge \quad \underline{\lim}\, \frac{\sin\tilde{\alpha}}{\sin\alpha} = \underline{\lim}\, \frac{\tilde{a}b\sin\tilde{\beta}}{\tilde{b}a\sin\beta} \\ &\ge \quad \underline{\lim}\, \frac{H(a)}{a}\frac{b}{H(b)} \cdot \underline{\lim}\, \frac{\sin\tilde{\beta}}{\sin\beta} \\ &\ge \quad \underline{\lim}\, \frac{\sin\tilde{\beta}}{\sin\beta} > \frac{1}{2}.\end{aligned}$$

□

We define $G_1(t) = \min\{t,\sqrt{t}\}$, $t \ge 0$, and recall definition (8) of the functions $F_k$.

LEMMA 4.3. *Consider a sequence $(\Delta_k)$ of Euclidean triangles, or of spherical triangles with circumscribed radius uniformly bounded by $R<\pi/2$. Let the angles of $\Delta_k$ be $\alpha_k \le \beta_k \le \gamma_k$, and assume that $\alpha_k \to 0$, $\beta_k \to \beta_0 \in (0,\pi)$, $\gamma_k \to \gamma_0 \in (0,\pi)$.*

*Put $H_k = G_1$ in the case of Euclidean triangles, and $H_k = F_k$ in the case of spherical triangles for $k \in \mathbb{N}$.*

*Then*

$$\underline{\lim}\, D(H_k,\Delta_k) > 1/2.$$

*Proof.* We use the notation of Lemmas 4.1 and 4.2 and their proofs. Note that the functions $H_k$, $k \in \mathbb{N}$, satisfy the assumptions of Lemma 4.2 in the Euclidean case, and the requirements 1–3 of Lemma 4.1 in the spherical case. So it is enough to prove

$$\underline{\lim}\, \tilde{\gamma}_k/\gamma_k > 1/2. \tag{53}$$

Our functions $H_k$ have the the following additional property

$$H_k(x)/H_k(y) \leq \max\{\sqrt{x/y}, \sqrt{\operatorname{chd} x/\operatorname{chd} y}\} \quad \text{for} \quad 0 < x \leq y. \tag{54}$$

Now we again drop the subscript $k$.

From (54), the asymptotic relations $a/b \to 0$ obtained in the proofs of Lemmas 4.1 and 4.2, and $\operatorname{chd} x \asymp x$ for $x > 0$ bounded away from $\pi$, we conclude $\tilde{a}/\tilde{b} = H(a)/H(b) \to 0$ both in the Euclidean and spherical case.

Then $\sin \tilde{\alpha} \to 0$ by the Law of Sines. Since $\tilde{\alpha}$ is the smallest angle in $\tilde{\Delta}$ this implies $\tilde{\alpha} \to 0$. Hence $\tilde{\beta} + \tilde{\gamma} \to \pi$. Since $\tilde{\gamma}$ is the largest angle in $\tilde{\Delta}$ this shows $\underline{\lim}\, \tilde{\gamma} \geq \pi/2$. Thus $\underline{\lim}\, \tilde{\gamma}/\gamma \geq \pi/(2\gamma_0)$ which gives (53). □

## 5. The case of two small angles. Euclidean triangles

Now we have to investigate what happens when two angles become small. As we are dealing with a spherical triangle whose circumscribed radius is bounded away from $\pi/2$, the presence of two small angles implies that the diameter of the triangle is small (see Lemma 3.5). So it is 'almost' Euclidean, which explains why we study Euclidean triangles here. The side distortion function is

$$G_k(t) = \min\{kt, \sqrt{t}\}, \quad t \geq 0, \quad k > 0. \tag{55}$$

As the angles of Euclidean triangles are scaling invariant it is enough to consider the case when $k = 1$ in (55). Notice the following three properties of $G = G_1$:

$$G(x + y) \leq G(x) + G(y), \quad x, y \geq 0. \tag{56}$$

This holds because $G$ is concave, increasing and $G(0) = 0$.

$$G^2(x + y) \geq G^2(x) + G^2(y), \quad x, y \geq 0. \tag{57}$$

The function $G^2$ is not convex, but nevertheless satisfies (57), which can be seen by consideration of four simple cases depending on the location of the three points $x, y$ and $x + y$ with respect to the point $1$.

$$G(x^2) = G^2(x), \quad x \geq 0. \tag{58}$$

In this section we consider Euclidean triangles $\Delta$ with sides $a, b, c$ and angles $\alpha, \beta, \gamma$ opposite to these sides. The transformed triangle $\tilde{\Delta} = G\Delta$ has sides $\tilde{a}, \tilde{b}, \tilde{c}$ and angles $\tilde{\alpha}, \tilde{\beta}, \tilde{\gamma}$ such that the sides or angles that correspond under the transformation are denoted by the same letter.

LEMMA 5.1. *If* $G: [0,\infty) \to [0,\infty)$ *is a concave increasing function satisfying* (56), (57), (58) *and* $G \not\equiv 0$, *then* $D(G, \Delta) > 1/2$ *for every Euclidean triangle* $\Delta$.

*Proof.* It follows from concavity that $G(t) > 0$ for $t > 0$. Then (57) implies that $G$ is strictly increasing.

We check the distortion of $\gamma$. By the Cosine of Half-Angles Formula [12, p. 60]

$$c = \sqrt{(a+b)^2 - 4ab\cos^2(\gamma/2)}. \tag{59}$$

We are going to prove the inequality $\tilde{\gamma} > \gamma/2$, which is equivalent to $\cos\tilde{\gamma} < \cos(\gamma/2)$. By the Law of Cosines this is equivalent to

$$\frac{G^2(a) + G^2(b) - G^2(c)}{2G(a)G(b)} < \cos\frac{\gamma}{2}. \tag{60}$$

Introducing $t = \cos(\gamma/2)$, $0 < t < 1$, we substitute (59) into (60) and use (58) to transform (60) to the following equivalent form

$$G^2(a) + G^2(b) - 2G(a)G(b)t < G\left((a+b)^2 - 4abt^2\right). \tag{61}$$

For fixed $a$ and $b$ the right-hand side of (61) is strictly concave (concave and not linear) with respect to $t$ as a composition of the strictly increasing concave function $G$ and the strictly concave function $t \mapsto (a+b)^2 - 4abt^2$. The left-hand side of (61) is linear in $t$. So by the Maximum Principle it is enough to verify the inequality on the boundary of the segment $[0,1]$. Moreover, it is enough to verify nonstrict inequality on the boundary, to conclude that (61) is a strict inequality.

For $t = 0$ we use (58) again and obtain $G^2(a) + G^2(b) \le G^2(a+b)$, which is the same as (57).

For $t = 1$ we also use (58) and obtain $(G(a) - G(b))^2 \le G^2(|a-b|)$ which is equivalent to (56). □

Again as in Section 3 the bound $1/2$ for the distortion in Lemma 4.1 is the best possible. The next lemma singles out the case when the distortion can be close to $1/2$.

LEMMA 5.2. *There exists an increasing function* $g: (0, 2\pi/3] \to (0,\infty)$, $0 < g(t) \le t$, $t \in (0, 2\pi/3]$, *with the following property. For every Euclidean triangle* $\Delta$ *whose smaller angles are* $\alpha$ *and* $\beta$ *we have for every* $k > 0$ $D(G_k, \Delta) \ge 1/2 + g(\alpha+\beta)$, *where* $G_k$ *is defined in* (55).

*Proof.* For every $s>0$ and every triangle $\Delta$ we denote by $s\Delta$ the triangle similar to $\Delta$ obtained by multiplying all sides of $\Delta$ by $s$.

As $G_k\Delta$ is similar to $G_1k^2\Delta$ we can assume without loss of generality that $k=1$.

It is enough to show that if the sum of the two smaller angles of a triangle $\Delta$ is bounded below by a constant $\mu\in(0,2\pi/3]$, then the distortion $D(G_1,\Delta)$ is bounded away from $1/2$.

Assuming this not the case, we find a sequence $(\Delta_k)$ of triangles with the sum of the two smaller angles

$$\alpha_k+\beta_k\geq\mu>0, \tag{62}$$

and

$$\lim D(G_1,\Delta_k)\leq 1/2. \tag{63}$$

Without loss of generality we have $\alpha_k\leq\beta_k$. Moreover, by selecting a subsequence we may assume that all three angles of $\Delta_k$ have limits as $k\to\infty$. Now we consider two cases.

*Case* 1. $\lim\alpha_k>0$. Then $\lim\beta_k>0$. We put $s_k:=(\mathrm{diam}\Delta_k)^{-1}$. By selecting a subsequence we may assume that $s_k$ tends to a limit, possibly infinite.

If $\lim s_k\in(0,\infty)$, then $\Delta_k\to\Delta$ for some (nondegenerate!) triangle $\Delta$, and from (63) we obtain $D(G_1,\Delta)\leq 1/2$, which contradicts Lemma 5.1.

If $\lim s_k=\infty$, then $s_k\Delta_k\to\Delta$ for some triangle $\Delta$. The triangles $G_1\Delta_k$ are similar to $s_kG_1\Delta_k$, so by (63)

$$D(G^0,\Delta)\leq 1/2, \tag{64}$$

where $G^0(t):=\lim_{k\to\infty}s_kG_1(s_k^{-1}t)=t$. So $G^0$ is a linear function, and it does not distort angles at all, which contradicts (64).

If $\lim s_k=0$, then again $s_k\Delta_k\to\Delta$ for some triangle $\Delta$. The triangles $G_1\Delta_k$ are similar to $s_k^{1/2}G_1\Delta_k$, so by (63)

$$D(G^\infty,\Delta)\leq 1/2, \tag{65}$$

where $G^\infty(t):=\lim_{k\to\infty}s_k^{1/2}G_1(s_k^{-1}t)=\sqrt{t}$. Obviously, $G^\infty$ is a monotone function which satisfies (56), (57) and (58). Then (65) contradicts Lemma 5.1. This completes Case 1.

*Case* 2. $\lim\alpha_k=0$. Then our assumption (62) implies that $\beta_k$ has a positive limit, which is strictly less than $\pi$, because $\beta_k$ is the second largest angle in $\Delta_k$. Moreover, it is clear that the limit of the largest angles $\gamma_k$ in $\Delta_k$ also belongs to $(0,\pi)$.

So we are in the situation of the Euclidean case of Lemma 4.3. Thus we obtain a contradiction to (63). ☐

Now we need a better estimate for the distortion of the largest angle in a triangle with two small angles.

LEMMA 5.3. *Let* $\Delta$ *be a Euclidean triangle with angles* $\alpha \leq \beta \leq \gamma$. *Let* $\tilde{\gamma}$ *be the angle in* $\tilde{\Delta} = G_k\Delta$, $k > 0$, *corresponding to the angle* $\gamma$ *in* $\Delta$, *where* $G_k$ *is defined in* (55). *Then*

$$\tilde{\gamma} \geq \pi/2 + O(\alpha+\beta)^2, \quad \alpha+\beta \to 0,$$

*uniformly with respect to* $\Delta$ *and* $k$.

*Proof.* As in Lemma 5.2 it is enough to consider the case $k = 1$. We have $a \leq b \leq c$. Put

$$\phi := \pi/2 - \tilde{\gamma}.$$

We want to estimate $\phi$ from above. From the Law of Sines

$$a = \frac{c\sin\alpha}{\sin(\alpha+\beta)} \quad \text{and} \quad b = \frac{c\sin\beta}{\sin(\alpha+\beta)}. \tag{66}$$

*Case* 1. $a \leq b \leq c \leq 1$. Since in this case $\tilde{\Delta} = \Delta$ and $\tilde{\gamma} = \gamma \to \pi$ as $\alpha + \beta \to 0$, the statement is obvious in this case.

*Case* 2. $a \leq b \leq 1 < c$. From the Law of Cosines and (66) we now obtain

$$\sin\phi = \cos\tilde{\gamma} = \frac{c(\sin^2\alpha + \sin^2\beta) - \sin^2(\alpha+\beta)}{2c\sin\alpha\sin\beta}. \tag{67}$$

The right-hand side of (67) is an increasing function of $c$ for fixed $\alpha$ and $\beta$, so we can use the estimate

$$c \leq \frac{\sin(\alpha+\beta)}{\sin\beta},$$

which follows from (66) and our assumption that $b \leq 1$. By substituting this estimate to (67) we obtain using $\alpha \leq \beta$

$$\begin{aligned}
\sin\phi &\leq \frac{\sin^2\alpha + \sin^2\beta - \sin(\alpha+\beta)\sin\beta}{2\sin\alpha\sin\beta} \\
&= \frac{\sin^2\alpha + \sin^2\beta - (\sin\alpha\cos\beta + \cos\alpha\sin\beta)\sin\beta}{2\sin\alpha\sin\beta} \\
&= \frac{\sin^2\alpha - \sin\alpha\sin\beta\cos\beta + O(\alpha^2\sin^2\beta)}{2\sin\alpha\sin\beta} \\
&= \frac{\sin\alpha}{2\sin\beta} - \frac{\cos\beta}{2} + O(\alpha\sin\beta) \\
&\leq O(\beta^2).
\end{aligned}$$

This completes Case 2.

*Case* 3. $a \le 1 < b \le c$. From the Law of Cosines and (66) we now obtain

$$\sin\phi = \cos\tilde{\gamma} = \frac{c\sin^2\alpha + \sin\beta\sin(\alpha+\beta) - \sin^2(\alpha+\beta)}{2\sin\alpha\sqrt{c\sin\beta\sin(\alpha+\beta)}}. \tag{68}$$

The right-hand side of (68) is an increasing function of $c$ for fixed $\alpha$ and $\beta$, so we can use the estimate

$$c \le \frac{\sin(\alpha+\beta)}{\sin\alpha}$$

which follows from (66) and our assumption $a \le 1$. By substituting this in the right-hand side of (68) we obtain,

$$\begin{aligned}
\sin\phi &\le \frac{\sin\alpha + \sin\beta - \sin(\alpha+\beta)}{2\sqrt{\sin\alpha\sin\beta}} \\
&= \frac{\sin\alpha + \sin\beta - \sin\alpha\cos\beta - \sin\beta\cos\alpha}{2\sqrt{\sin\alpha\sin\beta}} \\
&= \frac{\sin\alpha\sin^2(\beta/2) + \sin\beta\sin^2(\alpha/2)}{2\sqrt{\sin(\alpha/2)\sin(\beta/2)\cos(\alpha/2)\cos(\beta/2)}} \\
&= \frac{\sin(\alpha/2)\cos(\alpha/2)\sin^2(\beta/2) + \sin(\beta/2)\cos(\beta/2)\sin^2(\alpha/2)}{\sqrt{\sin(\alpha/2)\sin(\beta/2)\cos(\alpha/2)\cos(\beta/2)}} \\
&= \frac{\cos(\alpha/2)\sin^{3/2}(\beta/2)\sin^{1/2}(\alpha/2) + \cos(\beta/2)\sin^{3/2}(\alpha/2)\sin^{1/2}(\beta/2)}{\sqrt{\cos(\alpha/2)\cos(\beta/2)}} \\
&= O(\alpha+\beta)^2.
\end{aligned}$$

This completes Case 3.

*Case* 4. $1 < a \le b \le c$. In this case $G_1\Delta = G_\infty\Delta$, where $G_\infty$ was defined in (32). So we can use Lemma 3.2 which says that the smallest ratio $\tilde{\gamma}/\gamma$ will occur for the isosceles triangle of the same circumscribed radius as $\Delta$, and equal sides meeting at $\gamma$. This isosceles triangle has the same angle $\gamma$ and thus the same sum $\alpha+\beta$ as our triangle $\Delta$. So in this case we can assume without loss of generality that $\alpha = \beta$.

From the Law of Cosines and (66) we obtain

$$\begin{aligned}
\sin\phi &= \frac{\sin\alpha + \sin\beta - \sin(\alpha+\beta)}{2\sqrt{\sin\alpha\sin\beta}} \\
&= \frac{2\sin\beta - \sin(2\beta)}{2\sin\beta} = 2\sin^2(\beta/2) \\
&= O(\beta^2).
\end{aligned}$$

This completes Case 4 and the proof of Lemma 5.3. ☐

## 6. Conclusion of the proof of Theorem 1.5

*Proof of Theorem* 1.5. The proof is by contradiction. We assume that there exists $\varepsilon > 0$ such that for every $k \geq 1$ and for every $\delta > 0$ there exists a spherical triangle $\Delta$ with circumscribed radius

$$R(\Delta) \leq b_0 - \varepsilon \tag{69}$$

and

$$D(F_k, \Delta) \leq 1/2 + \delta.$$

In particular, putting $\delta = \delta_k = g^2(k^{-3})$ for $k \in \mathbb{N}$, where $g$ is the function from Lemma 5.2, there will be a triangle $\Delta = \Delta_k$ with the properties (69) and

$$D(F_k, \Delta_k) \leq 1/2 + g^2(k^{-3}). \tag{70}$$

We denote the sides of $\Delta_k$ by $a_k \leq b_k \leq c_k$, the angles of $\Delta_k$ by $\alpha_k \leq \beta_k \leq \gamma_k$, the corresponding angles of the transformed triangle $\Delta'_k := \operatorname{chd} \Delta_k$ by $\alpha'_k \leq \beta'_k \leq \gamma'_k$, and the angles of $\tilde{\Delta}_k = F_k \Delta_k = G_k \Delta'_k$ by $\tilde{\alpha}_k \leq \tilde{\beta}_k \leq \tilde{\gamma}_k$. Moreover, we denote the spherical circumscribed radius of $\Delta_k$ by $R_k = R(\Delta_k)$, and its spherical diameter by $d_k$. By Lemma 3.5 these quantities satisfy the following asymptotic relation

$$\alpha_k + \beta_k \asymp d_k / R_k, \quad k \to \infty. \tag{71}$$

It is also clear that

$$d_k \leq 2R_k. \tag{72}$$

By selecting a subsequence of $(\Delta_k)$, we may assume that all angles involved, as well as $(d_k)$, $(R_k)$ and $(d_k/R_k)$ have limits in $[0, \infty)$ for $k \to \infty$.

Now we consider the following cases 1–4. In all cases, except the very last one, we will show that

$$\underline{\lim}\, D(F_k, \Delta_k) > 1/2, \tag{73}$$

which will contradict (70) because $g(t) \to 0$ as $t \to 0$. Only in Case 4 do we use the full strength of (70) to obtain a contradiction.

*Case* 1. $\lim d_k > 0$.

Then $\lim R_k > 0$ by (72), and thus by (71)

$$\lim \alpha_k + \lim \beta_k =: m > 0. \tag{74}$$

If $\lim \alpha_k > 0$, the sequence $(\Delta_k)$ tends to a nondegenerate spherical triangle $\Delta$, and $F_k \Delta_k = F_\infty \Delta_k$ for $k$ large enough. So by Lemma 3.4 we have $\lim D(F_k, \Delta_k) = D(F_\infty, \Delta) > 1/2$ which gives (73).

If $\lim \alpha_k = 0$, then $\beta_0 := \lim \beta_k > 0$. Moreover, we claim that $\gamma_0 := \lim \gamma_k < \pi$. To prove the claim, we assume the contrary, that is, $\gamma_k \to \pi$. Then from Delambre's formula (sometimes credited to Gauss, [12, p. 162])

$$\cos(\gamma_k/2)\cos\Big((a_k - b_k)/2\Big) = \cos(c_k/2)\sin\Big((\alpha_k + \beta_k)/2\Big)$$

and the condition that the lengths of the sides of our triangles are bounded away from $\pi$ we obtain $\sin\Big((\alpha_k+\beta_k)/2\Big) \to 0$. By (74) this implies $\alpha_k + \beta_k \to 2\pi$. Then $\operatorname{area}(\Delta_k) = \alpha_k + \beta_k + \gamma_k - \pi \to 2\pi$. This is impossible, since the circumscribed radius $R_k$ of $\Delta_k$ is bounded away from $\pi/2$, the radius of a hemisphere.

So $0 < \beta_0 \leq \gamma_0 < \pi$, and (73) holds by the spherical case of Lemma 4.3.

*Case* 2. $\lim d_k = 0$ but $\lim\, d_k/R_k > 0$.
Then by (71) we again have (74). In addition $\lim R_k = 0$, so by Lemma 2.4 the distortion of angles coming from the projection $\Pi$ is negligible; that is, $D(\mathrm{chd}\,, \Delta_k) \to 1$. Now we apply Lemma 5.2, and (74) to conclude that

$$\underline{\lim}\, D(F_k, \Delta_k) \geq \underline{\lim}\, D(G_k, \Delta'_k) \geq 1/2 + g(m/2).$$

This proves (73).

*Case* 3. $\lim d_k = 0$, $\lim\, d_k/R_k = 0$, and $d_k/R_k = o(k^{-2})$ as $k \to \infty$.
In this case, we have $F_k\Delta_k = k(\mathrm{chd}\,\Delta_k)$ for large $k$. Moreover, $\alpha_k + \beta_k \to 0$ by (71). We claim that

$$\underline{\lim}\, D(F_k, \Delta_k) \geq 1. \tag{75}$$

Indeed, from Lemma 3.5 follows that $\alpha_k \leq \alpha'_k = \tilde{\alpha}_k$ and $\beta_k \leq \beta'_k = \tilde{\beta}_k$ for large $k$. Moreover, since the circumscribed radii of the triangles $\Delta_k$ are bounded away from $\pi/2$ we have $\alpha_k + \beta_k \asymp \alpha'_k + \beta'_k$ by Lemma 2.4. Thus, $\lim(\tilde{\alpha}_k + \tilde{\beta}_k) = \lim(\alpha'_k + \beta'_k) = \lim(\alpha_k + \beta_k) = 0$ which implies that $\tilde{\gamma}_k = \pi - \tilde{\alpha}_k - \tilde{\beta}_k \to \pi$. This proves (75), which is stronger than (73).

It remains to consider

*Case* 4. $\lim d_k = 0$, $\lim d_k/R_k = 0$, and $\underline{\lim}\, k^2 d_k/R_k > 0$.
In this case, by (71) we have that for some $m > 0$

$$\underline{\lim}\, k^2(\alpha_k + \beta_k) =: 7m > 0, \tag{76}$$

and $\alpha_k + \beta_k \to 0$. First we consider the distortion of the smaller angles $\alpha_k$ and $\beta_k$. By Lemma 3.5 we have $\alpha_k \leq \alpha'_k$ and $\beta_k \leq \beta'_k$ for large $k$. Now we apply Lemma 5.2 to $\Delta'_k$ and obtain for large $k$

$$\begin{aligned} \tilde{\alpha}_k/\alpha_k &\geq \tilde{\alpha}_k/\alpha'_k \geq 1/2 + g(\alpha'_k + \beta'_k) \geq 1/2 + g(\alpha_k + \beta_k) \\ &\geq 1/2 + g(6mk^{-2}) \geq 1/2 + g(k^{-3}), \end{aligned} \tag{77}$$

and similarly

$$\tilde{\beta}_k/\beta_k \geq 1/2 + g(k^{-3}). \tag{78}$$

Now we estimate the distortion of the largest angle $\gamma_k$. By the Spherical Law of Cosines for Angles [12, p. 153], we have

$$\cos\gamma_k = -\cos\alpha_k\cos\beta_k + \sin\alpha_k\sin\beta_k\cos d_k.$$

Together with the estimate $d_k = O(\alpha_k+\beta_k)$, which follows from (71) this gives

$$\begin{aligned} 2\sin^2\left(\frac{\pi-\gamma_k}{2}\right) &= 2\cos^2(\gamma_k/2) = 1+\cos\gamma_k \\ &= 1-\cos\alpha_k\cos\beta_k+\sin\alpha_k\sin\beta_k\cos d_k \\ &= (\alpha_k+\beta_k)^2/2 + O(\alpha_k+\beta_k)^4, \end{aligned}$$

which implies

$$\gamma_k = \pi - (\alpha_k+\beta_k) + O(\alpha_k+\beta_k)^2. \tag{79}$$

On the other hand, using Lemma 5.3 and $\alpha_k \asymp \alpha_k'$, $\beta_k \asymp \beta_k'$ we obtain

$$\tilde{\gamma}_k \geq \pi/2 + O(\alpha_k'+\beta_k')^2 = \pi/2 + O(\alpha_k+\beta_k)^2.$$

Comparing this with (79), and taking into account our assumption (76) we get for large $k$

$$\tilde{\gamma}_k/\gamma_k \geq 1/2 + (\alpha_k+\beta_k)/(2\pi) + O(\alpha_k+\beta_k)^2 \geq 1/2 + mk^{-2} \geq 1/2 + k^{-3}.$$

Together with (77), (78) and $g(t) \leq t$, this implies $D(F_k,\Delta_k) \geq 1/2 + g(k^{-3})$ for large $k$ which contradicts (70) because $g(x)\to 0$ as $x\to 0$. □

## 7. Proof of Theorem 1.6

We will be using Lemmas 2.4, 3.5, 4.1 and 4.2. In addition we need the following modifications of Lemmas 4.3, 5.1, 5.2 and 5.3.

We recall that our side distortion function in Theorem 1.6 is

$$F_k^*(t) = \min\{k\,\mathrm{chd}\,t, 1\}.$$

If for $k>0$ we define

$$G_k^*(t) = \min\{kt, 1\}, \quad t \geq 0,$$

then we have $F_k^* = G_k^* \circ \mathrm{chd}$.

First we prove an analog of Lemma 4.3.

Lemma 7.1. *Consider a sequence* $(\Delta_k)$ *of Euclidean triangles, or of spherical triangles with circumscribed radius uniformly bounded by* $R < \pi/2$.

*Let the angles of* $\Delta_k$ *be* $\alpha_k \le \beta_k \le \gamma_k$, *and assume that* $\alpha_k \to 0$, $\beta_k \to \beta_0 \in (0, \pi)$, $\gamma_k \to \gamma_0 \in (0, \pi)$.

*Put* $H_k = G_1^*$ *in the case of Euclidean triangles, and* $H_k = F_k^*$ *in the case of spherical triangles for* $k \in \mathbb{N}$.

*Then*

$$\underline{\lim}\, D(H_k, \Delta_k) > 1/3.$$

*Proof.* We use the notation of Lemmas 4.1 and 4.2 and their proofs. Note that the functions $H_k$, $k \in \mathbb{N}$, satisfy the assumptions of Lemma 4.2 in the Euclidean case, and the requirements 1–3 of Lemma 4.1 in the spherical case. So the lower limits of the ratios $\tilde{\alpha}_k/\alpha_k$ and $\tilde{\beta}_k/\beta_k$ are at least $1/2$. Since $\tilde{\gamma}_k$ is the largest angle in the Euclidean triangle $\tilde{\Delta}_k$ we must have $\tilde{\gamma}_k \ge \pi/3$. Hence $\underline{\lim}\, \tilde{\gamma}_k/\gamma_k \ge \pi/(3\gamma_0) > 1/3$. □

Now we have to prove the results, similar to those of Section 5, for the function $G_1^*$.

LEMMA 7.2. $D(G_1^*, \Delta) > 1/3$ *for every Euclidean triangle* $\Delta$.

*Proof.* Denote the sides of $\Delta$ by $a \le b \le c$. Then the the largest side of $\tilde{\Delta} = G_1^* \Delta$ is $G_1^*(c)$ and the largest angle of $\tilde{\Delta}$ is $\tilde{\gamma}$. So $\tilde{\gamma} \ge \pi/3 > \gamma/3$. It remains to check the distortion of the smaller angles $\alpha$ and $\beta$. We consider the following four cases.

*Case* 1. $a \le b \le c \le 1$. This case is trivial because $G_1^* \Delta$ is similar to $\Delta$.

*Case* 2. $a, b \le 1 < c$. We don't use the convention $a \le b$ in this case. We claim that $\tilde{\alpha} \ge \alpha$ and $\tilde{\beta} \ge \beta$. By symmetry it is enough to check this for $\alpha$, and show that $\tilde{\alpha} \ge \alpha$, which is equivalent to $\cos \tilde{\alpha} \le \cos \alpha$. By the Law of Cosines this is equivalent to

$$\frac{G_1^*(c)^2 + G_1^*(b)^2 - G_1^*(a)^2}{2G_1^*(c) G_1^*(b)} \le \frac{c^2 + b^2 - a^2}{2bc},$$

which is the same as

$$\frac{1 + b^2 - a^2}{2b} \le \frac{c^2 + b^2 - a^2}{2bc},$$

or

$$c + cb^2 - ca^2 \le c^2 + b^2 - a^2, \quad \text{or} \quad (c-1)b^2 \le (c-1)c + (c-1)a^2,$$

or

$$b^2 \le c + a^2.$$

This is true since $b^2 \le 1 < c$ in this case.

*Case* 3. $a \leq 1 < b \leq c$. In this case $\tilde{\Delta}$ is isosceles with two sides equal to 1, and $\tilde{\beta} = \tilde{\gamma}$. Thus $\tilde{\beta} \geq \pi/3 > \beta/3$. It remains to check the distortion of $\alpha$. We verify that $\tilde{\alpha} \geq \alpha$, which is equivalent to $\cos\tilde{\alpha} \leq \cos\alpha$, or

$$\frac{1+1-a^2}{2} \leq \frac{b^2+c^2-a^2}{2bc}.$$

The last inequality simplifies to $2bc + a^2 \leq b^2 + c^2 + bca^2$ which is true since $2bc \leq b^2 + c^2$ and $bc \geq 1$.

*Case* 4. $1 \leq a \leq b \leq c$. This case is trivial because $\tilde{\Delta}$ is equilateral. □

The following is a counterpart of Lemma 5.2.

LEMMA 7.3. *There exists an increasing function* $g^*\colon (0, 2\pi/3] \to (0, \infty)$, $0 < g^*(t) \leq t$, $t \in (0, 2\pi/3]$, *with the following property. For every Euclidean triangle* $\Delta$ *whose smaller angles are* $\alpha$ *and* $\beta$ *we have* $D(G_k^*, \Delta) \geq (1/3) + g^*(\alpha + \beta)$ *for every* $k > 0$.

*Proof.* This is derived from Lemmas 7.2 and 7.1, Euclidean case, in the same way as Lemma 5.2 is derived from Lemmas 5.1 and 4.3, Euclidean case. Only minor modifications are necessary. Namely, $G_1$ should be replaced by $G_1^*$ and the constant 1/2 by 1/3 throughout the proof. In the last paragraph of Case 1 the function $G^\infty$ is defined by $G^\infty(t) = \lim G_1^*(s_k^{-1}t) \equiv 1$; so $G^\infty\tilde{\Delta}$ is equilateral and we get a contradiction to the modified version of (65).

To deal with Case 2 we use the Euclidean case of Lemma 7.1 to obtain a contradiction. □

Now we need the following trivial counterpart of Lemma 5.3.

LEMMA 7.4. *If* $\gamma$ *is the largest angle in a Euclidean triangle* $\Delta$, *then for the corresponding angle* $\tilde{\gamma}$ *in* $\tilde{\Delta} := G_k^*\Delta$ *we have* $\tilde{\gamma} \geq \pi/3$.

*Proof.* This is clear since $\tilde{\gamma}$ is the largest angle in $\tilde{\Delta}$. □

*Proof of Theorem* 1.6. The proof is by contradiction. We assume that there exists $\varepsilon > 0$ such that for every $k \geq 1$ and for every $\delta > 0$ there exists a spherical triangle $\Delta$ with circumscribed radius

$$R(\Delta) \leq \pi/2 - \varepsilon \tag{80}$$

and

$$D(F_k^*, \Delta) \leq 1/3 + \delta.$$

In particular, putting $\delta = \delta_k = (g^*(k^{-3}))^2$ for $k \in \mathbb{N}$, where $g^*$ is the function from Lemma 7.3, there will be a triangle $\Delta = \Delta_k$ with the properties (80) and

$$D(F_k^*, \Delta_k) \leq 1/3 + \Big(g^*(k^{-3})\Big)^2. \tag{81}$$

We denote the angles of $\Delta_k$ by $\alpha_k \le \beta_k \le \gamma_k$, the corresponding angles of the transformed triangle $\Delta'_k := \operatorname{chd} \Delta_k$ by $\alpha'_k \le \beta'_k \le \gamma'_k$, and of $\tilde{\Delta}_k = F^*_k \Delta_k = G^*_k \Delta'_k$ by $\tilde{\alpha}_k \le \tilde{\beta}_k \le \tilde{\gamma}_k$. Moreover, we denote the spherical circumscribed radius of $\Delta_k$ by $R_k = R(\Delta_k)$, and its spherical diameter by $d_k$. By Lemma 3.5 we have

$$\alpha_k + \beta_k \asymp d_k/R_k, \quad k \to \infty, \tag{82}$$

and the following holds

$$d_k \le 2R_k. \tag{83}$$

By selecting a subsequence of $(\Delta_k)$, we may assume that all angles involved, as well as $(d_k)$, $(R_k)$ and $(d_k/R_k)$ have limits in $[0, \infty)$ for $k \to \infty$.

Now we consider the following cases 1–4. In all cases, except the very last one we will show that

$$\underline{\lim}\, D(F^*_k, \Delta_k) > 1/3, \tag{84}$$

which will contradict (81) because $g^*(t) \to 0$ as $t \to 0$. Only in Case 4 do we use the full strength of (81) to obtain a contradiction.

*Case* 1. $\lim d_k > 0$. Then $\lim R_k > 0$ by (83), and thus by (82)

$$\lim \alpha_k + \lim \beta_k =: m > 0. \tag{85}$$

If $\lim \alpha_k > 0$, the sequence $(\Delta_k)$ tends to a nondegenerate spherical triangle $\Delta$, and $F^*_k \equiv 1$ on the sides of $\Delta_k$ for $k$ large enough. This means that the triangles $\tilde{\Delta}_k$ are equilateral, and all their angles are equal to $\pi/3$. Thus (84) follows.

If $\lim \alpha_k = 0$, then $\beta_0 := \lim \beta_k > 0$. As in Case 1 of the proof of Theorem 1.5 we see that $\gamma_0 := \lim \gamma_k < \pi$. Hence $0 < \beta_0 \le \gamma_0 < \pi$. So (84) holds by Lemma 7.1, spherical case.

*Case* 2. $\lim d_k = 0$ but $\lim d_k/R_k > 0$.
Then by (82) we again have (85). In addition $\lim R_k = 0$, so by Lemma 2.4 the distortion from the projection $\Pi$ is negligible; that is, $D(\operatorname{chd}, \Delta_k) \to 1$. Now we apply Lemma 7.3 and (85) to conclude that

$$\underline{\lim}\, D(F^*_k, \Delta_k) \ge \underline{\lim}\, D(G^*_k, \Delta'_k) \ge 1/3 + g^*(m/2).$$

This proves (84).

*Case* 3. $\lim d_k = 0$, $\lim d_k/R_k = 0$, and $d_k/R_k = o(k^{-2})$ as $k \to \infty$.
In this case we have

$$\underline{\lim}\, D(F^*_k, \Delta_k) \ge 1,$$

which is proved in the same way as in Case 3 of the proof of Theorem 1.5.

It remains to consider

*Case* 4. $\lim d_k = 0$, $\lim d_k/R_k = 0$, and $\underline{\lim}\, k^2 d_k/R_k > 0$.

In this case, by (82) we have that for some positive $m > 0$

$$\underline{\lim}\, k^2(\alpha_k + \beta_k) =: 10m > 0. \tag{86}$$

First we consider the distortion of the smaller angles $\alpha_k$ and $\beta_k$. By Lemma 3.5 we have $\alpha_k \leq \alpha'_k$ and $\beta_k \leq \beta'_k$ for large $k$. Now we apply Lemma 7.3 to $\Delta'_k$ and obtain for large $k$

$$\begin{aligned} \tilde{\alpha}_k/\alpha_k &\geq \tilde{\alpha}_k/\alpha'_k \geq 1/3 + g^*(\alpha'_k + \beta'_k) \geq 1/3 + g^*(\alpha_k + \beta_k) \\ &\geq 1/3 + g^*(9mk^{-2}) \geq 1/3 + g^*(k^{-3}), \end{aligned} \tag{87}$$

and similarly

$$\tilde{\beta}_k/\beta_k \geq 1/3 + g^*(k^{-3}). \tag{88}$$

Now we estimate the distortion of the largest angle $\gamma_k$. Again we have (79) from Section 6:

$$\gamma_k = \pi - (\alpha_k + \beta_k) + O(\alpha_k + \beta_k)^2. \tag{89}$$

On the other hand, Lemma 7.4 gives

$$\tilde{\gamma}_k \geq \pi/3.$$

Comparing this with (89), and taking into account our assumption (86) we get for large $k$

$$\tilde{\gamma}_k/\gamma_k \geq 1/3 + (\alpha_k + \beta_k)/(3\pi) + O(\alpha_k + \beta_k)^2 \geq 1/3 + mk^{-2} \geq 1/3 + k^{-3}.$$

Together with (87), (88) and $g^*(t) \leq t$ this implies $D(F^*_k, \Delta_k) \geq 1/3 + g^*(k^{-3})$ which contradicts (81). □

Technische Universität Braunschweig, Germany
*E-mail address*: m.bonk@tu-bs.de

Purdue University, West Lafayette, IN
*E-mail address*: eremenko@math.purdue.edu